\documentclass[12pt]{article}
\usepackage{e-jc+}

\title{Classification of Polynomial Mappings Between Commutative Groups}
\author{Uwe Schauz\\%\thanks{Thanks to the editors of this wonderful journal!}\\
\small Department of Mathematical Science\\[-0.8ex]
\small Xi’an Jiaotong-Liverpool University\\[-0.8ex]
\small Suzhou 215123, China\\[-0.8ex]
\small \texttt{uwe.schauz@xjtlu.edu.cn}}

\date{\printday\\%\dateline{\_\_\_ \_\_, \_\_\_\_}{\llap{--------------\ }%
\small Mathematics Subject Classifications: 05E99, 13F20, 41A05}

\addtolength\marginparsep{-2pt}

\usepackage{fixltx2e}
\usepackage[ansinew]{inputenc}%ß als Eingabe statt "s u.s.w.
\usepackage{datetime}%definiert \xxivtime, \ampmtime, ... .(Muss vor ngerman kommen)
\usepackage{ngerman}%neue Trennregeln.
\selectlanguage{USenglish}
\nonfrenchspacing%geringfügig größere Leerzeichenabstände der Satzzeichen.

\usepackage{eurosym}%DIN-Euro \euro \EUR1.20
\usepackage{textcomp}%TS1-Zeichencodierung z.B. \textdollar $.
\usepackage{url}
\usepackage[only,varolessthan,varogreaterthan,olessthan,ogreaterthan,shortleftarrow,shortrightarrow]{stmaryrd}

\usepackage{epic}
\usepackage{eepic}
\usepackage[%draft,%dvips,
debugshow]{graphicx}
\usepackage{amsmath}
\usepackage{amsthm}
\usepackage{extpfeil}
\usepackage{amssymb}%z.B. $\Box$ ,enthält die latexsym-Zeichen.
\usepackage{stmaryrd}%z.B. $\llbracket$
\usepackage{bm}%more general and robust implementation of \boldsymbol.

\usepackage{enumerate}
\usepackage{afterpage}
\usepackage{xspace}
\setlength\mathsurround{.2em}
\emergencystretch.03\textwidth%zusätzlicher Leerraum bei schwierigen Umbrüchen 1em="M"

\newtheoremstyle{Theorem}{.7\baselineskip}{1\baselineskip}{\itshape}{}{\bfseries}{.}{ }{}
\theoremstyle{Theorem}
\newtheorem{Satz}{Theorem}[section]

\newtheorem{Korollar}[Satz]{Corollary}

\newtheorem{Lemma}[Satz]{Lemma}

\newtheoremstyle{Definition}{.6\baselineskip}{.8\baselineskip}{}{}{\bfseries}{.}{ }{}
\theoremstyle{Definition}

\theoremstyle{definition}
\newtheorem{Beispiel}[Satz]{Example}

\theoremstyle{remark}
\newtheorem{Bemerkung}[Satz]{Remark}

\newenvironment{sequation}{\begin{small}\begin{equation}}{\end{equation}\end{small}}
\newenvironment{sequation*}{\begin{small}\begin{equation*}}{\end{equation*}\end{small}}
\newenvironment{Beweis}[1][Proof]{\begin{proof}[#1]}{\end{proof}}
\newcommand\ps{\small}
\newenvironment{proofsize}[1]{\begin{small}#1}{\end{small}}
\newcommand\psize[1]{\begin{proofsize}#1\end{proofsize}}

\newcommand\textps[1]{\text{\psize{#1}}}

\newcommand\´{\kern 1pt}%0.08em=1.44mu
\newcommand\ms{\hspace{\mathsurround}}
\newcommand\noms{\hspace{-\mathsurround}}
\newcommand\vsp{\vspace{1ex}}

\newcommand\printday{\today\xspace}
\newcommand\Rand[1]{%\mbox{}
  \marginpar{\raggedleft\scriptsize\hspace{0pt}#1}}%
\renewcommand{\(}{\noms$}
\renewcommand{\)}{\noms$}
\renewcommand\frac[2]{\genfrac{}{}{.4pt}{}{#1}{#2}}%EJC erzeugt zu breite Linien.
\renewcommand\tfrac[2]{\genfrac{}{}{.4pt}{1}{#1}{#2}}
\renewcommand\dfrac[2]{\genfrac{}{}{.4pt}{0}{#1}{#2}}
\renewcommand\atop[2]{\genfrac{}{}{0pt}{}{#1}{#2}}

\newcommand\mathRand[1]{\hspace{\mathsurround}\Rand{#1}\nolinebreak\noms}
\def\rand #1"#2"{\mathRand{\(#2\)}#1#2}
\def\randd #1"#2"#3\randd#4"#5"{\mathRand{\(#2\), \(#5\)}#1#2#3#4#5}

\newcommand\eqby[2][=]%
  {\ensuremath{\overset{\makebox[0pt]{\ensuremath{\smash[t]{\scriptstyle#2}}}}{#1}}}
\newcommand\F[1][\ ]{\mathbb{F}_{\!#1}}

\newcommand\Q{\mathbb{Q}}
\newcommand\Z{\mathbb{Z}}
\newcommand\N{\mathbb{N}}

\renewcommand\sb{\subseteq}%vorher: XX\sb{index}=XX_{index}
\newcommand\downmto{\rotatebox[y=7pt]{-90}{\ensuremath{\mto}}}
\newcommand\mto{\mapsto}
\newcommand\lmto{\longmapsto}
\DeclareMathOperator*\lto{\longrightarrow}
\newcommand\To{\Rightarrow}
\newcommand\lTo{\Longrightarrow}

\newcommand\lEqi{\Longleftrightarrow}
\newcommand\nach{\mathbin{\scalebox{.55}{\raisebox{.218em}{\ensuremath{\pmb\varolessthan}}}}}
\newcommand\fa{\forall\,}

\newcommand\DP{\colon\discretionary{\!\kern -.17em}{}{}}
\newcommand\mitsymbol{\textup{\textbrokenbar}}
\renewcommand\mit{\,\ \discretionary{\mitsymbol}{}{}\mitsymbol\ \,}
\renewcommand\div{\mathrel{\bigm\lfloor\!\!\!\bigm\lfloor}}
\newcommand\vid{\mathrel{\bigm\rfloor\!\!\!\bigm\rfloor}}
\newcommand\ndiv{\mathrel{\;\!\div\hspace{-12pt}\kern0pt\lower2pt%
  \hbox{\ensuremath{^\diagup}}\!}}
\newcommand\ndivps{\mathrel{\;\!\div\hspace{-9pt}\kern0pt\lower2pt%
  \hbox{\ensuremath{^\diagup}}\!}}
\newcommand\nvid{\mathrel{\;\!\vid\hspace{-12pt}\kern0pt\lower2pt%
  \hbox{\ensuremath{^\diagup}}\!}}
\newcommand\nvidps{\mathrel{\;\!\vid\hspace{-9pt}\kern0pt\lower2pt%
  \hbox{\ensuremath{^\diagup}}\!}}
\newcommand\Div{\mathrel{\Bigm\lfloor\hspace{-6pt}\Bigm\lfloor}}

\newcommand\nDiv{\mathrel{\,\Div\hspace{-12.13pt}\diagup}}

\newcommand\n[1][n]{{\phantom{|}}^{\!\!\!\!#1}}

\providecommand\abs[1]{\lvert#1\rvert}

\newcommand\Spnb[1]%
{\bigl\langle\hspace{-5pt}\bigl\langle#1\bigr\rangle\hspace{-5pt}\bigr\rangle}

\DeclareMathOperator\Id{Id}
\DeclareMathOperator\End{End}

\newcommand\x{\chi}

\renewcommand\b{\beta}
\newcommand\ä{\alpha}
\newcommand\tx{\skew{2}{\tilde}{x}}
\newcommand\hx{\skew{2}{\hat}{x}}
\renewcommand\d{\delta}

\newcommand\Chi{\kern0pt\lower-2.5pt\hbox{\ensuremath\chi\!}}

\newcommand\cop[1]{\lfloor#1\rfloor}%{[#1]}%
\begingroup\makeatletter\ifx\SetFigFont\undefined
\def\x#1#2#3#4#5#6#7\relax{\def\x{#1#2#3#4#5#6}}%
\expandafter\x\fmtname xxxxxx\relax \def\y{splain}%
\ifx\x\y   % LaTeX or SliTeX?
\gdef\SetFigFont#1#2#3{%
  \ifnum #1<17\tiny\else \ifnum #1<20\small\else
  \ifnum #1<24\normalsize\else \ifnum #1<29\large\else
  \ifnum #1<34\Large\else \ifnum #1<41\LARGE\else
     \huge\fi\fi\fi\fi\fi\fi
  \csname #3\endcsname}%
\else
\gdef\SetFigFont#1#2#3{\begingroup
  \count@#1\relax \ifnum 25<\count@\count@25\fi
  \def\x{\endgroup\@setsize\SetFigFont{#2pt}}%
  \expandafter\x
    \csname \romannumeral\the\count@ pt\expandafter\endcsname
    \csname @\romannumeral\the\count@ pt\endcsname
  \csname #3\endcsname}%
\fi
\fi\makeatother\endgroup

\begin{document}
\maketitle

\begin{abstract}
Some polynomials $P$ with rational coefficients give rise to well defined maps between
cyclic groups, $\Z_q\lto\Z_r$, $x+q\Z\lmto P(x)+r\Z$. More generally, there are polynomials
in several variables with tuples of rational numbers as coefficients that induce maps between
commutative groups. We characterize the polynomials with this property, and classify all
maps between two given finite commutative groups that arise in this way. We also provide
interpolation formulas and a Taylor-type theorem for the calculation of polynomials that
describe given maps.
\end{abstract}

\section{Introduction}\label{sec.intr}
The polynomial
\begin{equation}
P\,:=\,-\tfrac{1}{8}X^4+\tfrac{3}{4}X^3-\tfrac{7}{8}X^2-\tfrac{3}{4}X+1
\end{equation}
induces a map $\Z_3\lto \Z_9$ in a canonical way. At first, it gives rise to a map $\Z\lto \Q$
which actually is integer valued, as one can show. Second, this map $\Z\lto \Z$ induces the
map $\Z\lto \Z_9$, $x\lmto P(x)+9\Z$. Finally, it turns out that our new map is even
\(3\)"~periodic,
\begin{equation}
P(x+3)\,\equiv\,P(x)\pmod{9}\,\ .
\end{equation}
So, we obtain a well defined map
\begin{equation}\label{eq.PMod}
P\DP\Z_3\lto \Z_9\ , \quad\ x=\hx+3\Z\lmto P(x):=P(\hx)+9\Z\,\ .
\end{equation}
Moreover, the map $P$ is a kind of \emph{Lagrange Polynomial}. For
$x\in\Z_3:=\Z/3\Z$\mathRand{\(\Z_q\)}
\begin{equation}
P(x)\ =\ \dbinom{0}{x}_{\!\!3,9}\,:=\,\begin{cases}
 1+9\Z&\text{if $\,x=0$\,,}\\
 0+9\Z&\text{if $\,x\neq0$\,.}
\end{cases}
\end{equation}

As any function is a linear combination of Lagrange Functions, we see that any map
$\Z_3\lto \Z_9$ can be represented by a polynomial in $\Q[X]$ of degree at most $4$. Of
course, not every polynomial over $\Q$ gives rise to a well defined map $\Z_3\lto \Z_9$, but
we have enough polynomials to obtain all maps. This is not true for maps $\Z_q\lto\Z_r$ in
general, it holds only if $q$ and $r$ are powers of a common prime $p$. In contrast, if $q$ is
coprime to $r$ then only constant maps can be described by a rational polynomial. These
two extremal results are the antagonistic forces that determine the general case. To reduce
the question of representability of maps between $\Z_q$ and $\Z_r$, with general $q$ and
$r$, to these two cases, we split the domain $\Z_q$ and the codomain $\Z_r$ into cyclic
\(p\)"~groups, and decompose polynomial maps $\Z_q\lto\Z_r$ into maps between the cyclic
\(p\)-factors. More precisely, we write the domain $\Z_q$ as direct product of $n$ cyclic
\(p\)"~groups, and introduce one variable $X_j$ for each of them. Similarly, the codomain
$\Z_r$ is split into $t$ cyclic \(p\)"~groups, and the coefficients of our polynomials are split
into rational \(t\)"~tuples accordingly. These decompositions can be done without changing
the property of polynomial representability (Theorem\,\ref{sz.split}), but they bring us closer
to our two extremal cases. Decompositions also allow the treatment of maps between
arbitrary finite commutative groups $A$ and $B$. Our main result, Theorem\,\ref{sz.per}\´,
says that, if $A_1,A_2,\dotsc,A_t$ and $B_1,B_2,\dotsc,B_t$ are the (possibly trivial or
noncyclic) primary components of $A$ and $B$ (corresponding to the different prime
divisors $p_1,p_2,\dotsc,p_t$ of $\abs{A}\abs{B}$), then the polynomial representable maps
in $B^A$ are exactly the maps
\begin{equation}
(f_1,\dotsc,f_t)\in B_1^{A_1}\!\times\dotsm\times B_t^{A_t}\ ,\ \quad
(a_1,\dotsc,a_t)\mto(f_1(a_1),\dotsc,f_t(a_t))\,\ .
\end{equation}

Partial results in this direction, but sometimes over more general rings and sometimes in the
language of iterated differences, were obtained in several other papers, e.g.\ in
\cite{ccf,eimu,fr1,fr2,hr,is,ji,jpsz,la,lei,mm,must,wi}, and in the books \cite{cc,na}\´. Our
paper differs from most of these investigations in that we do not restrict ourselves to the
case where domain and codomain coincide. We also do not just count or determine the
isomorphy type of modules of certain maps.

One important point in our definitions of polynomial maps on residue classes $x+q\Z\in\Z_q$
is the independence from representatives $x$. It guarantees that composition of polynomial
maps corresponds to substitution (based on Corollary\,\ref{cor.comp}). For example, if
$Q,P\in\Q[X]$ describe maps $\Z_{27}\lto \Z_3$ and $\Z_3\lto \Z_9$, respectively, then
$P(Q(X))\in\Q[X]$ describes the composed map $\Z_{27}\lto \Z_9$. In particular, there is a
polynomial of degree at most $\deg(P)\deg(Q)$ that describes the composed map. In fact,
the degree is the important new parameter that comes with the polynomials. It plays the key
role in important theorems like the Combinatorial Nullstellensatz about which we want to
write another paper. One frequently used tool in this connection are Lagrange Polynomials
$1-X^{q-1}$ over a finite fields $\F[q]$, with the important property to turn nonzeros into
zeros, and vice versa, $1-x^{q-1}=0\lEqi x\neq0$. For example, they are used in the usual
proofs of Chevalley and Warning's Theorem, as e.g.\ in \cite[Corollary\,3.5]{schAlg}\´. Our
Lagrange Polynomials in Theorem\,\ref{sz.lagr}\´, together with the degree restriction in
Theorem\,\ref{sz.md}\´, will provide here a straightforward way to generalize such results.
More classical applications of Lagrange Polynomials can be found in interpolation theory, as
e.g.\ in the proof of our Theorem\,\ref{sz.ip}\´.

In this paper, we talk about \emph{polyfracts}, which are nothing else than integer valued
polynomials, but in the notation of binomial polynomials with modular coefficients. These
notions are introduced in Section\,\ref{sec.int}\´, together with a Taylor-type theorem based
on iterated differences. In this section the domain is always $\Z$ or $\Z^n$. Only in
Section\,\ref{sec.per}\´, we start with periodicity and with cyclic and commutative groups as
domain. This section also contains our main results and, at the very end, an example that
demonstrates our classification theorems.

% -----------------------------------------------------------------------------
\section{Integer- and Group-Valued Polynomials}\label{sec.int}

In this section, we introduce group"=valued polynomials in $n$ variables on the domain
$\Z^n$\!. We start with polynomials in one variable with values in a cyclic group.

\subsection{Polyfractal Maps Into Cyclic Groups}

It is well known that the integer valued polynomials $P\in\Q[X]$ (\´i.e.\ $P(z)\in\Z$ for all
$z\in\Z$) are precisely the \emph{binomial polynomials} over $\Z$, or
\(\Z\)"~\emph{polyfracts}, as we say here (see Theorem\,\ref{sz.ur}). These are \(\Z\)"~linear
combinations of \emph{binomial monomials}, or \emph{(mono)fracts}; which are, for
$\d\in\N$, defined by
 \rand\begin{equation}
"\dbinom{X}{\d}"\,:=\,\dfrac{X(X-1)\dotsm(X-\d+1)}{\d!}\quad\text{and}\quad\dbinom{X}{0}\,:=\,1\,\
.
\end{equation}
One example of a polyfract would be the polynomial $P$ from the introduction,
\begin{equation}
P\,:=\ -\tfrac{1}{8}X^4+\tfrac{3}{4}X^3-\tfrac{7}{8}X^2-\tfrac{3}{4}X+1
\ =\ -3\dbinom{X}{4}+\dbinom{X}{2}-\dbinom{X}{1}+\dbinom{X}{0}\,\ .
\end{equation}
This \emph{monofractal expansion}, and the well known fact that binomial coefficients are
natural numbers, explains why $P$ is integer valued. The  monofractal representability also
suggests the notation \rand\begin{equation} "\Z\tbinom{X}{\Z}"\,:=\,\{\,P\in\Q[X]\mit
P|_\Z\in\Z^\Z\,\}
\end{equation}
for the subring of polyfracts. Since, in this paper, we mostly view the function values modulo
$r$\!, we also may view the coefficients of monofracts modulo $r$\!. We write
\rand$"\Z_r\tbinom{X}{\Z}"$ for the set of such \(\Z_r\)"~linear combinations of monofracts,
where we usually allow the extremal cases
  \rand\rand\begin{equation}
"\Z_1":=\{0\}\quad\text{and}\quad"\Z_0":=\Z\,\ .
\end{equation}
Actually, the factor ring
\begin{equation}
\Z_r\tbinom{X}{\Z}\,=\,\Z\tbinom{X}{\Z}\!\bigm/\!r\Z\tbinom{X}{\Z}
\end{equation}
is the best possible choice here, as the following theorem shows\´:

\begin{Satz}\label{sz.1}
For $\geq0$, the ring homomorphism
$$
\Z\tbinom{X}{\Z}\lto\Z_r^\Z\ ,\quad Q\lmto Q|_\Z+r\Z^\Z
$$
has kernal $r\Z\binom{X}{\Z}$ and induces the monomorphism
$$
\Z_r\tbinom{X}{\Z}\lto\Z_r^\Z\ ,\quad\bar Q\lmto\bar Q|_\Z\,\ .
$$
Hence, to any map $f\in\Z_r^\Z$ there exists at most one polyfract $P\in\Z_r\tbinom{X}{\Z}$
that describes (interpolates) this map, $P|_\Z=f$. We do not have to make a difference
between polyfracts $P$ and polyfractal maps $P|_\Z$, and may view polyfracts as special
maps,
$$
\Z_r\tbinom{X}{\Z}\,\sb\,\Z_r^\Z\ .
$$
\end{Satz}

\begin{Beweis}
Let $Q\in\Z\tbinom{X}{\Z}$. Obviously, if $Q\in r\Z\tbinom{X}{\Z}$ then $Q$ is mapped to
zero. To prove the converse, assume $Q$ is mapped to zero,
\begin{sequation}
Q|_\Z+r\Z^\Z\,=\,0\,\in\,\Z^\Z\!\bigm/\!r\Z^\Z\,=\,\Z_r^\Z\ ,
\end{sequation}%
so that
\begin{sequation}
Q|_\Z\,\in\,r\Z^\Z\,.
\end{sequation}%
Then $\tfrac{1}{r}Q$ %as polynomial in $\Q[X]$
is integer valued, and $\tfrac{1}{r}Q|_\Z$ can be described by a polyfract
$R\in\Z\tbinom{X}{\Z}$,
\begin{sequation}
R|_\Z\,:=\,\tfrac{1}{r}Q|_\Z\,\in\,\Z^\Z\,.
\end{sequation}%
From this follows
\begin{sequation}
Q\,=\,rR\,\in\,r\Z\tbinom{X}{\Z}\,\ ,
\end{sequation}%
as the coefficients of polyfracts $Q$ over $\Q$ are uniquely determined by the infinitely
many function values $Q(x)$, $x\in\Z$, like those of polynomials. This proves the main
statement about the kernal. As
$\Z_r\tbinom{X}{\Z}=\Z\tbinom{X}{\Z}\!\bigm/\!r\Z\tbinom{X}{\Z}$, we also obtain the
described embedding of $\Z_r\tbinom{X}{\Z}$ in $\Z_r^\Z$.
\end{Beweis}

Surprisingly, as we have seen, our modulo $r$ treatment of coefficients allows us to identify
a polyfract $P$ with its polyfractal map $P|_{\Z}$. We are used to the opposite from the
study of polynomials, where a statement like this would be completely wrong. For example,
the two polynomials $X^2-X$ and $0$ in $\Z_2[X]$ both describe the zero map
$0\DP\Z\lto\Z_2$. There is not just one unique polynomial representative. One may wonder
about this, as any polynomial can be written as polyfract, and both, $X^2-X$ and $0$, should
give rise to a polyfract describing the zero map. However, viewing the two different
polynomials as polyfracts modulo $2$ makes them equal,
\begin{equation}
0\,=\,0\dbinom{X}{2}\,=\,2\dbinom{X}{2}\,=\,X^2-X\ \quad\text{in $\Z_2\tbinom{X}{\Z}$.}
\end{equation}
In other words,
\begin{equation}
X^2-X\,\in\,2\Z\tbinom{X}{\Z_2}\ \ \quad\text{and}\quad\
\Z_2\tbinom{X}{\Z}\,\cong\,\Z\tbinom{X}{\Z}\!\bigm/\!2\Z\tbinom{X}{\Z}\,\ .
\end{equation}
We will see that it can be convenient to identify $P\in\Z_r\binom{X}{\Z}$ with $P|_{\Z}$.
Based on the identification, we always have two points of view. For example, if we view two
polyfracts $P$ and $Q$ as maps, then addition and multiplication is just pointwise,
\begin{equation}
(P+Q)(x)\,=\,P(x)+Q(x)\quad\text{and}\quad(PQ)(x)\,=\,P(x)Q(x)\,\ .
\end{equation}
If we view $P$ and $Q$ as linear combinations of monofracts, then addition is
coefficientwise, and equally simple as before. However, multiplication in the polyfractal
representation is quite complicated. One may produce multiplication tables for monofracts
and use distributivity, or one can follow the following five steps\´:

\begin{enumerate}
\item Replace the coefficients of $P$ and $Q$ by integer representatives to obtain
    representatives of $P$ and $Q$ in $\Z\tbinom{X}{\Z}$.
\item Expand them as polynomials in $\Q[X]$.
\item Perform the multiplication in $\Q[X]$.
\item Express the result $PQ\in\Q[X]$ as \(\Z\)"~polyfract,
    $PQ=c_m\tbinom{X}{m}+\dotsb+c_1\tbinom{X}{1}+c_0\tbinom{X}{0}$, $m:=\deg(PQ)$.
    To calculate $c_m$ divide $PQ$ by $\tbinom{X}{m}\in\Q[X]$ inside $\Q[X]$. The
    reminder yields the other coefficients recursively.
\item Reduce the coefficients $c_j$ modulo $r$.
\end{enumerate}\medskip

\subsection{Polyfractal Maps Into Groups}\label{sec.mi}

More generally than above, we may allow the elements of any finitely generated
\(\Z\)"~module (i.e.\ commutative group) $B$ as coefficients, and write
\rand$"B\tbinom{X}{\Z}"$ for the set of \(B\)"~polyfracts. Since $B$ is the direct product of
cyclic groups,
\begin{equation}
B\,=\,\Z_{r_1}\times\Z_{r_2}\times\dotsb\times\Z_{r_t}\qquad\text{with numbers $r_i\geq0$,}
\end{equation}
our set $B\tbinom{X}{\Z}$ can be written as
\begin{equation}
(\Z_{r_1}\!\times\Z_{r_2}\!\times\dotsb\times\Z_{r_t})\dbinom{X}{\Z}
\,\,=\,\,\Z_{r_1}\dbinom{X}{\Z}\times\Z_{r_2}\dbinom{X}{\Z}\times\dotsb\times
\Z_{r_t}\dbinom{X}{\Z}\,\ .
\end{equation}
Here, we always identify the both sides, e.g.\ we identify
\begin{equation}
(2,3)\dbinom{X}{2}+(4,0)\dbinom{X}{1}+(1,5)\,=\,\Bigl(2\tbinom{X}{2}+4\tbinom{X}{1}
+1\,\,,\,\,3\tbinom{X}{2}+5\Bigr)\,\ .
\end{equation}
From this representation $B\tbinom{X}{\Z}$ inherits a ring structure, since the factors
$\Z_{r_i}\tbinom{X}{\Z}$ carry a ring structure already. We just have to add and multiply
componentwise, viewing a polyfract over
$\Z_{r_1}\!\times\Z_{r_2}\!\times\dotsb\times\Z_{r_t}$ as \(n\)"~tuples of polyfracts over
$\Z_{r_1},\Z_{r_2},\dotsc,\Z_{r_t}$, respectively. For any \(\Z\)"~module homomorphism
\begin{equation}
\phi\DP C\lto B\,\ ,
\end{equation}
we also denote with $\phi$ the \(\Z\)"~module homomorphisms
\begin{equation}
\phi\DP C\tbinom{X}{\Z}\lto B\tbinom{X}{\Z}
\ \quad\text{and}\quad\phi\DP C^\Z\lto B^\Z\,\ ,
\end{equation}
defined by coefficientwise, respectively pointwise, application of $\phi\DP C\lto B$. The
definitions of these two induced maps go well together, if we identify polyfracts with
polyfractal maps, as suggested in Theorem\,\ref{sz.1}\´. In other words, the diagram
\begin{equation}
\begin{split}
 P\!&\ \xmapsto{\qquad\qquad\quad\ \´\´\ }\phi(P)\\
 &\!\!\!\!\smash{\downmto}\qquad
  \,\ \raisebox{3.5pt}[0em][0em]{\ensuremath{\circlearrowright}}\qquad\qquad\quad\,
  \smash{\downmto}\\[-5pt]
 P|_\Z\!&\ \lmto\ \phi\nach P|_\Z=\phi(P)|_\Z\\
\end{split}
\end{equation}
commutes. However, in this paper $\phi$ will usually just be the invers to the \emph{splitting
isomorphism} (\emph{Chinese Reminder Isomorphism})
\begin{equation}
\begin{split}
\Z_r\,\lto&\ \Z_{r_1}\times\Z_{r_2}\times\dotsb\times\Z_{r_t}\\
x+r\Z\ \lmto&\,(x+r_1\Z\´,x+r_2\Z\´,\dotsc,x+r_t\Z)\,\ ,
\end{split}
\end{equation}
corresponding to the factorization
\begin{equation}
r\,=\,r_1\,r_2\,\dotsm\,r_t
\end{equation}
of $r\geq2$ into pairwise coprime prime powers, and to the decomposition
\begin{equation}
\Z_r\,=\,\tfrac{r}{r_1}\Z_r\oplus\dotsb\oplus\tfrac{r}{r_t}\Z_r
\end{equation}
of $\Z_r$ into Sylow subgroups. Explicitly,
\begin{equation}
\phi\DP x_1{+}r_1\Z\´,\dotsc,x_t{+}r_t\Z\,\,\lmto\,\,x_1s_1\dfrac{r}{r_1}+\dotsb
+x_ts_t\dfrac{r}{r_t}\,+\,r\Z
\end{equation}
where, $s_1,\dotsc,s_t$ are numbers with
\begin{equation}
s_1\dfrac{r}{r_1}+\dotsb+s_t\dfrac{r}{r_t}
\ =\ \gcd\Bigl(\dfrac{r}{r_1},\dotsc,\dfrac{r}{r_t}\Bigr)\ =\ 1\,\ .
\end{equation}

Originally, we thought that polyfractal maps into cyclic groups $\Z_r$ are the right object of
study, but then it turned out that prime powers play a special rule in our theory. Therefore, it
will be convenient to break cyclic codomains $\Z_r$ down into cyclic \(p\)"~groups. This is
what we need splitting isomorphisms for. However, this prime power factorization works for
all finite commutative groups $B=\Z_{r_1}\!\times\Z_{r_2}\!\times\dotsb\times\Z_{r_t}$ as
codomain. Since splitting isomorphisms may also be applied to their factors $\Z_{r_i}$, we
can always go over to the case where $r_1,r_2,\dotsc,r_t$ are already prime powers. We
just have the drop the unnecessary assumption that our prime powers $r_1,r_2,\dotsc,r_t$
are pairwise coprime, which would hold for cyclic $B$.

\subsection{Polyfractal Maps in Several Variables}

We can introduce multivariate polyfracts $P=P(X_1,X_2,\dotsc,X_n)$ over finitely generated
commutative groups $B=\Z_{r_1}\!\times\Z_{r_2}\!\times\dotsb\times\Z_{r_t}$ in exactly the
same way as multivariate polynomials are introduced, as sums of products of polyfracts in
one variable, together with a distributive law. This means, they are finite linear combinations
of monofracts in $n$ variables, which are just products of the form
\begin{equation}
\dbinom{X_1,X_2,\dotsc,X_n}{\d_1,\d_2,\dotsc,\d_n}
\,:=\,\dbinom{X_1}{\d_1}\dbinom{X_2}{\d_2}\dotsm\dbinom{X_n}{\d_n}\,\ ,
\end{equation}
with
\begin{equation}
\binom{X_1,X_2,\dotsc,X_n}{\d_1,\d_2,\dotsc,\d_n}=1\quad\text{if $n=0$.}
\end{equation}
Hence, our polyfracts $P=P(X_1,X_2,\dotsc,X_n)$ can be written as
\begin{equation}
P\,=\,\sum_{\d\in\N^n}P_\d\dbinom{X_1,X_2,\dotsc,X_n}{\d_1,\d_2,\dotsc,\d_n}
\end{equation}
with only finitely many nonvanishing coefficients $P_\d\in B$. We denote with
\begin{equation}\label{eq.B}
\begin{split}
B\dbinom{X_1,X_2,\dotsc,X_n}{\Z^n}\,&=\,(\Z_{r_1}\!\times\dotsb\times
\Z_{r_t})\dbinom{X_1,X_2,\dotsc,X_n}{\Z^n}\\
\,&=\,\Z_{r_1}\dbinom{X_1,X_2,\dotsc,X_n}{\Z^n}\times\dotsb\times
\Z_{r_t}\dbinom{X_1,X_2,\dotsc,X_n}{\Z^n}
\end{split}
\end{equation}
the ring of all such polyfracts. If $r_1=r_2=\dotsb=r_t=0$, this is precisely the subring of all
polynomials in $\Q^t[X_1,X_2,\dotsc,X_n]$ that map any $(x_1,x_2,\dotsc,x_n)\in\Z^n$ into
$\Z^t$\!, which follows immediately from the following generalization of
\cite[Theorem\,2.1]{na}\´:

\begin{Satz}\label{sz.ur}
The polyfracts $P\in\Z\tbinom{X_1,X_2,\dotsc,X_n}{\Z^n}$ are precisely the integer valued
polynomials in $\Q[X_1,X_2,\dotsc,X_n]$.
\end{Satz}

\begin{Beweis}
Due to the combinatorial interpretation of binomial coefficients $\tbinom{d}{\d}$ as number of
\(\d\)"~subsets of a \(d\)"~set, it is obvious that polyfracts
$P\in\Z\tbinom{X_1,X_2,\dotsc,X_n}{\Z^n}$ are integer valued. To prove the converse, let
$P\in\Q[X_1,X_2,\dotsc,X_n]$ be integer valued. We may assume $P\neq0$, and that the
theorem was already proven for all polynomials $P'$\! of lower total degree or of same total
degree but with viewer monomials of that degree. Now, let $X_1^{\d_1}X_2^{\d_2}\dotsm
X_n^{\d_n}$ be such a monomial of maximal total degree in $P$\!,
\begin{sequation}
\d_1+\d_2+\dotsb+\d_n\,=\,\deg(P)\,\ .
\end{sequation}%
We will see in Equation\,\eqref{eq.do} in the next section that the coefficient $P_{[\d]}\neq0$
of this monomial can be written as
\begin{sequation}
P_{[\d]}\,:=\,\dfrac{P_\d}{\d!}\,\ ,
\end{sequation}%
with a $P_\d\in\Z$ and
\begin{sequation}
\d!:=\d_1!\,\d_2!\dotsm\d_n!\,\ .
\end{sequation}%
Hence, the polynomial
\begin{sequation}
P'\,:=\,P-P_\d\dbinom{X_1,X_2,\dotsc,X_n}{\d_1,\d_2,\dotsc,\d_n}
\end{sequation}%
will not contain $X_1^{\d_1}X_2^{\d_2}\dotsm X_n^{\d_n}$ any more. It will contain viewer
monomials of degree $\deg(P)$, and none of higher degree. Hence, by our induction
hypothesis,
\begin{sequation}
P'\,\in\,\Z\binom{X_1,X_2,\dotsc,X_n}{\Z^n}\,\ ,
\end{sequation}%
and
\begin{sequation}
P\,=\,P'+P_\d\dbinom{X_1,X_2,\dotsc,X_n}{\d_1,\d_2,\dotsc,\d_n}
\,\in\,\Z\dbinom{X_1,X_2,\dotsc,X_n}{\Z^n}
\end{sequation}%
follows.
\end{Beweis}

\begin{Korollar}\label{cor.comp}
Let $Q\in\Z\tbinom{X}{\Z}$ and $P\in\Z\tbinom{X_1,X_2,\dotsc,X_n}{\Z^n}$, then the
composed polynomial $Q(P)$ is a polyfract in $\Z\tbinom{X_1,X_2,\dotsc,X_n}{\Z^n}$ and
$\deg(Q(P))=\deg(Q)\deg(P)$.
\end{Korollar}

\begin{Beweis}
The composed polynomial in $\Q[X_1,X_2,\dotsc,X_n]$ is integer valued, i.e.\ polyfractal.
The equation $\deg(Q(P))=\deg(Q)\deg(P)$ holds, as it holds for polynomials.
\end{Beweis}

Theorem\,\ref{sz.ur} shows that the situation in several variables is precisely as it is for
polyfracts in just one variable. This is why Theorem\,\ref{sz.1} generalizes without difficulty to
the following result, which we formulated over arbitrary commutative groups here\´:

\begin{Satz}
Let $B$ be a finite commutative group. The map
$$
B\tbinom{X_1,X_2,\dotsc X_n}{\Z^n}\lto B^{\Z^n}\ ,\quad P\lmto P|_{\Z^n}
$$
is \emph{injective}. To any map $f\in B^{\Z^n}$ there exists at most one polyfract $P\in
B\tbinom{X_1,X_2,\dotsc X_n}{\Z^n}$ that describes (interpolates) this map, $P|_{\Z^n}=f$.
We do not have to make a difference between polyfracts $P$ and polyfractal maps
$P|_{\Z^n}$, and may view polyfracts as special maps,
$$
B\tbinom{X_1,X_2,\dotsc X_n}{\Z^n}\,\sb\,B^{\Z^n}\ .
$$
\end{Satz}

We also provide the following interesting little theorem\´:

\begin{Satz}\label{sz.fI}
Let $B$ be a finitely generated commutative group and
$$
\small{P=\sum_{\d\in\N^n}P_\d\dbinom{X_1,X_2,\dotsc,X_n}{\d_1,\d_2,\dotsc,\d_n}
\,\in\, B\dbinom{X_1,X_2,\dotsc,X_n}{\Z^n}}
$$
For grids $[d]:=\{0,1,\dotsc,d_1\}\times\dotsm\times\{0,1,\dotsc,d_n\}$ the following are
equivalent\´:
\begin{enumerate}[(i)]
\item $(\d\mto P_\d)|_{[d]}\,\equiv\,0$.
\item $(x\mto P(x))|_{[d]}\,\equiv\,0$.
\end{enumerate}
\end{Satz}

\begin{Beweis}
The proof is based on the simple fact that
\begin{sequation}
\dbinom{x}{\d}=0\quad\text{if $0\leq x<\d$.}
\end{sequation}%
The implication $(i)\To(ii)$ follows immediately from this. To prove $\neg(i)\To\neg(ii)$, let
$\d\in [d]$ be minimal with $P_\d\neq0$, i.e.
\begin{sequation}
P_{\d'}\neq0\ \ \lTo\ \ \d'\geq\d\,\ ,
\end{sequation}%
then
\begin{sequation}
P(\d)\,=\,\sum_{\d'\geq\d}P_{\d'}\dbinom{\d_1,\d_2,\dotsc,\d_n}{\d_1',\d_2',\dotsc,\d_n'}
\,=\,P_\d\,\neq\,0\,\ .
\end{sequation}%
\end{Beweis}

\subsection{The Discrete Derivative and Taylor's Theorem}

Our polyfractal notation also goes well together with an alternative type of derivative. For
polynomials $f\in\Q[X]$, and also for functions $f$ from a cyclic additive group into an
additive group, we define the \emph{(discrete) derivative} or \emph{difference function}
\rand$"\Delta f"$ via
\begin{equation}
\Delta f(x)\,:=\,f(x+1)-f(x)\,\ .
\end{equation}
If we have $n$ variables then we also have $n$ \emph{difference operators}
$\underset{^1}\Delta,\underset{^2}\Delta,\dotsc,\underset{^n}\Delta$. Applied to polynomials
\begin{equation}
P(X_1,X_2,\dotsc,X_n)\,=\,\sum_{\d\in\N^n}P_{[\d]}\,X_1^{\d_1}X_2^{\d_2}\dotsm X_n^{\d_n}\,\ ,
\end{equation}
these operations yield
\begin{equation}\label{eq.do}
\d!\,P_{[\d]}\,=\,\underset{^1}\Delta^{\!\d_1}\underset{^2}\Delta^{\!\d_2}\dotsm
\underset{^n}\Delta^{\!\d_n}P(0)
\end{equation}
for all $\d$ with
\begin{equation}
\d_1+\d_2+\dotsb+\d_n\,\geq\,\deg(P)\,\ .
\end{equation}
This follows readily from
\begin{equation}
\Delta X^\ell\,=\,\ell X^{\ell-1}+\,\textps{\emph{lower degrees}}\,\ ,
\end{equation}
where the terms of lower degree make it difficult to derive an equation for the coefficients
$P_{[\d]}$ with $\d_1+\d_2+\dotsb+\d_n<\deg(P)$. However, life is much easier if we work
with polyfractal representations. From Pascal's rule follows, for $k,\ell\in\N$,
\begin{equation}\label{eq.Dmi}
\Delta\dbinom{X}{\ell+1}\,=\,\dbinom{X}{\ell}\ \ ,\qquad
\Delta\dbinom{X}{0}\,=\,0
\end{equation}
and
\begin{equation}\label{eq.dmi}
\Bigl.\Delta^{\!k}\dbinom{X}{\ell}\Bigr|_{X=0}\,=\,\begin{cases}
 1&\text{if $\,k=\ell$\,,}\\
 0&\text{if $\,k\neq\ell$\,,}
\end{cases}
\end{equation}
which implies a kind of Taylor Theorem (as in \cite[Equation\,(\(3.7\))]{hr})\´:

\begin{Satz}\label{sz.tay}
Let $r\in\N$ and $f\DP\Z\lto\Z_r$ be a function with $\Delta^{\!d+1}f\equiv0%\pmod{r}
$, then
$$
f(x)\,\equiv\,\sum_{\d=0}^d\Delta^{\!\d} f(0)\dbinom{x}{\d}\,\ .%\pmod{r}\,\ .
$$
\end{Satz}

\begin{Beweis}
Both sides of the equation have the same \hbox{$(d{+}1)^{\textup{st}}$} derivative, by
Equation\,\eqref{eq.Dmi}. Furthermore, all lower derivatives coincide in $0$, by
Equation\,\eqref{eq.dmi}. Therefore, it suffices to prove that two functions $f$ and $g$ are
equal if
\begin{sequation}
\Delta f=\Delta g\ \quad\text{and}\quad f(0)=g(0)\,\ .
\end{sequation}%
This can be shown step by step as follows\´:
\begin{alignat}{3}
f(1)&\,\ =\,\,\ &\Delta f(0)+f(0)&\,\ =\,\,\ &\Delta g(0)+g(0)&\,\ =\,\ g(1)\,\ ,\notag\\
f(2)&\,\ =\,\,\ &\Delta f(1)+f(1)&\,\ =\,\,\ &\Delta g(1)+g(1)&\,\ =\,\ g(2)\,\ ,\\
\vdots\quad&&\vdots\qquad\ &&\vdots\qquad\ &\quad\qquad\vdots\notag
\end{alignat}\\[-2.2em]
and
\begin{alignat}{3}
f(-1)&\,\ =\,\ &f(0)-\Delta f(-1)&\,\ =\,\,\ &g(0)-\Delta g(-1)&\,\ =\,\ g(-1)\,\ ,\notag\\
f(-2)&\,\ =\,\ &f(-1)-\Delta f(-2)&\,\ =\,\,\ &g(-1)-\Delta g(-2)&\,\ =\,\ g(-2)\,\ ,\\
\vdots\quad&&\vdots\qquad\qquad&&\vdots\qquad\qquad&\quad\qquad\vdots\qquad\,.
\notag
\end{alignat}\\[-2em]
\end{Beweis}

From this theorem we see that a function $f\DP\Z\lto\Z_r$ is a polyfract, i.e.\ an integer
valued polynomial, of degree at most $d$, if and only if $\Delta^{\!d+1}f\equiv0$. One may
take this as definition for polynomials, as in many other papers, e.g.\ in \cite{hr,la,lei,must}\´,
see also Theorem\,\ref{sz.Equ}. Obviously, we also have the following \(n\)"~dimensional
version of this theorem\´:

\begin{Satz}\label{sz.tay2}
Let $r\in\N$ and $f\DP\Z\lto\Z_r$ be a map with
$\underset{^j}\Delta^{\!d_j+1}f\,\equiv\,0$ for $j=1,2,\dotsc,n$, then
$$
f(x)\,\equiv\,\sum_{\d_1=0}^{d_1}\,\sum_{\d_2=0}^{d_2}\dotsb\sum_{\d_n=0}^{d_n}
\underset{^1}\Delta^{\!\d_1}\underset{^2}\Delta^{\!\d_2}\dotsm
\underset{^n}\Delta^{\!\d_n}f(0)\,
\dbinom{x_1,x_2,\dotsc,x_n}{\d_1,\d_2,\dotsc,\d_n}\,\ .
$$
\end{Satz}

Our Taylor theorem may also be applied to polyfracts
\begin{equation}
P\,:=\,P_0\tbinom{X}{0}+P_1\tbinom{X}{1}+\dotsb+P_d\tbinom{X}{d}\,\ ,
\end{equation}
yielding
\begin{equation}
P_\d\,=\,\Delta^{\!\d} P(0)\,=\,\tbinom{\d}{0}P(\d)-\tbinom{\d}{1}P(\d-1)+
\dotsb+(-1)^\d\tbinom{\d}{\d}P(0)
\end{equation}
for all $\d\in\N$. From this we obtain another simple but remarkable insight, supplementing
our Theorem\,\ref{sz.fI}\´:

\begin{Bemerkung}\label{rem.inf}
The first $q$ function values $P(0),P(1),\dotsc,P(q-1)$ of a polyfract $P$ in one variable
suffice to calculate the first $q$ coefficients $P_0,P_1,\dotsc,P_{q-1}$, and vice versa. This
is of particular interest with respect to  \(q\)"~\emph{periodic} polyfracts $P$\!, as defined
and studied in the next section. In this case, all information about the function $x\lmto P(x)$
is already contained in the first $q$ evaluations of $P$\!, i.e.\ also in its first $q$ coefficients.
Hence, we may call the $q$ first coefficients the \emph{information coefficients}. The other
\emph{periodicity coefficients} are just required to make the polyfractal map periodic, and
are uniquely determined by the information coefficients.
\end{Bemerkung}

This remark gave a preliminary insight into periodic polyfracts, which we study more
systematically in the following section.

\section{Polynomials Between Groups}
\label{sec.per}

In this section, we move from the infinite domain $\Z^n$ to finite domains
$\Z_{q_1}\!\times\Z_{q_2}\!\times\dotsm\times\Z_{q_n}$, i.e., to arbitrary finite commutative
groups. Throughout, $p$ denotes a prime number.

\subsection{Periodicity and Polyfractal Maps On Groups}

First, we study \(q\)"~periodic functions $f\DP\Z\lto\Z_r:=\Z/r\Z$ (with $\Z_0=\Z$), i.e.
\begin{equation}
f(x+q)\,=\,f(x)\ \quad\text{for all $x\in\Z$.}
\end{equation}
In other words, our map $f\DP\Z\lto\Z_r$ is \(q\)"~periodic if and only if
\mathRand{\(\Delta_q, T\)}
\begin{sequation} \Delta_qf\,\equiv\,0\quad\text{where}\ \
 \Delta_q\ :=\ T^q-\Id\quad\text{with}\ \ T\DP f\lmto T f\ ,\,\ T f(x)\,:=\,f(x+1)\,\ .
\end{sequation}%
Our \(q\)"~periodic functions are basically the same as maps from $\Z_q$ to $\Z_r$. If a
\,\(q\)"~periodic function $f\DP\Z\lto\Z_r$ is given, then a corresponding map
$f\DP\Z_q\lto\Z_r$, again denoted $f$, is well defined by
\begin{equation}
f(x+q\Z)\,:=\,f(x)\,\ .
\end{equation}
This reinterpretation works also the other way around. With this slight interpretational shift
towards \(q\)"~periodicity in mind, we may view $\Z_r^{\Z_q}$ as subset of $\Z_r^{\Z}$. We
have
\begin{equation}
\Z_r\,=\,\Z_r^{\Z_1}\,\sb\,\Z_r^{\Z_q}\,\sb\,\Z_r^{\Z_{mq}}
\,\sb\,\Z_r^{\Z}\qquad\text{for $m=1,2,\dotsc$.}
\end{equation}
Also worth mentioning is that, if $r=q$, any polynomial $P\in\Z_q[X]$ gives rise to a map
$P|_{\Z_q}\DP\Z_q\lto\Z_q$, so that the corresponding polynomial map
$P|_\Z\DP\Z\lto\Z_q$ is \(q\)"~periodic. However, this is not true for polyfracts over $\Z_q$
(with $q>0$), e.g.
\begin{equation}
P(0+q)\,\neq\,P(0)\,\in\,\Z_q\ \quad\text{if}\quad P(X):=\dbinom{X}{q}\,\in\Z_q\tbinom{X}{\Z}\,\ .
\end{equation}
We write
 \rand\begin{equation}
"\Z_r\dbinom{X}{\Z_q}"\,:=\,\Bigl\{P\in\Z_r\tbinom{X}{\Z}\!\mit\!P\textps{\ is
\(q\)"~periodic}\Bigr\}\,=\,\Z_r\dbinom{X}{\Z}\cap\´\Z_r^{\Z_q}
\end{equation}
for the set of \(q\)"~periodic polyfracts over $\Z_r$, i.e.\ those \(\Z_r\)"~linear combinations
of monofracts that happen to be \(q\)"~periodic. Polyfracts $P$ in $n$ variables
$X_1,X_2,\dotsc,X_n$ may be \((q_1,q_2,\dotsc,q_n)\)"~periodic, i.e.\ \(q_j\)"~periodic with
respect to $X_j$, $j=1,2,\dotsc,n$. We denote the set of \((q_1,q_2,\dotsc,q_n)\)"~periodic
polyfracts over $\Z_r$ by\mathRand{\(\Z_r\!\dbinom{X_1,\dotsc}{\Z_{q}^n\,}\)}
\begin{equation}
\Z_r\dbinom{X_1,\,X_2,\,\dotsc,\,X_n}{\Z_{q_1}\!\times\Z_{q_2}\!\times\dotsm\times\Z_{q_n}}
\,:=\,\Z_r\dbinom{X_1,\dotsc,X_n}{\Z^n}\cap\´\Z_r^{\Z_{q_1}\!\times\Z_{q_2}
 \!\times\dotsm\times\Z_{q_n}}\,\ .
\end{equation}
Here, as before, we may replace the cyclic group $\Z_r$ in these definitions with any finitely
generated commutative group $B$. The cyclic case is just our basic core. We always can go
over to the general case by taking cartesian products, either coefficient"=wise on the level of
polyfracts, as in Equation\,\eqref{eq.B}, or on the level of maps using that
\begin{equation}
\Z_{r_1}^{\Z_{q_1}\!\times\dotsm\times\Z_{q_n}}\times\dotsm\times\Z_{r_t}^{\Z_{q_1}
 \!\times\dotsm\times\Z_{q_n}}
\,=\,(\´\Z_{r_1}\!\times\dotsm\times\Z_{r_t})^{\Z_{q_1}\!\times\dotsm\times\Z_{q_n}},
\end{equation}
where $t$ maps $f_i\DP\Z_{q_1}\!\times\Z_{q_2}\!\times\dotsm\times\Z_{q_n}\lto\Z_{r_i}$,
$i=1,2,\dotsc,t$, are combined to one map
\begin{equation}
\begin{split}
(f_1,f_2,\dotsc,f_t)\,\DP\,\Z_{q_1}\!\times\Z_{q_2}\!\times\dotsm\times\Z_{q_n}
\,\,\lto&\ \,\Z_{r_1}\!\times\Z_{r_2}\!\times\dotsm\times\Z_{r_t}\\
\mathbf{x}\,\ \lmto&\ \bigl(f_1(\mathbf{x}),f_2(\mathbf{x}),\dotsc,f_t(\mathbf{x})\bigr)\,\ .
\end{split}
\end{equation}

Based on our Taylor Theorems on can easily prove the following equivalence. It is a typical
result, a bit more general, in some aspects, than similar results in other papers, as e.g.\
\cite{la}\´:

\begin{Satz}\label{sz.Equ}
Let $d_1,d_2,\dotsc,d_n,q_1,q_2,\dotsc,q_n\geq0$ and $r_1,r_2,\dotsc,r_t\geq2$. For
maps $f\DP\Z_{q_1}\!\times\Z_{q_2}\!\times\dotsm\times\Z_{q_n}\lto\,
\,\Z_{r_1}\!\times\Z_{r_2}\!\times\dotsm\times\Z_{r_t}$ the following are equivalent\´:
\begin{enumerate}[(i)]
\item \(\underset{^j}\Delta^{\!d_j+1}f\,\equiv\,0\ms\) for $j=1,2,\dotsc,n$\,.
\item \(f\)\ms arises from a polynomial $P\in\Q^t[X_1,X_2,\dotsc,X_n]$ with restricted
    partial degrees $\deg_jP\leq d_j$, $j=1,2,\dotsc,n$\,.
\item \(f\)\ms arises from a polyfract
    $P\in(\Z_{r_1}\!\times\Z_{r_2}\!\times\dotsm\times\Z_{r_t})\dbinom{X_1,X_2,\dotsc,X_n}{\Z^n}$
    with restricted partial degrees $\deg_jP\leq d_j$, $j=1,2,\dotsc,n$\,.
\end{enumerate}
\end{Satz}

%\begin{Satz}\label{sz.Equ}
%Let $\d_1,\d_2,\dotsc,\d_n,q_1,q_2,\dotsc,q_n\geq0$ and $r_1,r_2,\dotsc,r_t\geq2$. For
%maps $f\DP\Z_{q_1}\!\times\Z_{q_2}\!\times\dotsm\times\Z_{q_n}\lto\,
%\,\Z_{r_1}\!\times\Z_{r_2}\!\times\dotsm\times\Z_{r_t}$ the following are equivalent\´:
%\begin{enumerate}[(i)]
%\item \(\underset{^1}\Delta^{\!\d_1}\underset{^2}\Delta^{\!\d_2}\dotsm\underset{^n}
%    \Delta^{\!\d_n}f\,\equiv\,0\)\ms.
%\item \(f\)\ms arises from a polynomial in $\Q^t[X_1,X_2,\dotsc,X_n]$ that does not
%    contain monomials $X_1^{\e_1}X_2^{\e_2}\dotsm X_n^{\e_n}$ with $\e_1\geq\d_1$,
%    $\e_2\geq\d_2$, \dots, $\e_n\geq\d_n$\,.
%\item \(f\)\ms arises from a polyfract in
%    $(\Z_{r_1}\!\times\Z_{r_2}\!\times\dotsm\times\Z_{r_t})
%    \dbinom{X_1,X_2,\dotsc,X_n}{\Z^n}$ that does not contain monofracts
%    $\dbinom{X_1,X_2,\dotsc,X_n}{\e_1,\,\e_2,\,\dotsc,\,\e_n}$ with $\e_1\geq\d_1$,
%    $\e_2\geq\d_2$, \dots, $\e_n\geq\d_n$\,.
%\end{enumerate}
%\end{Satz}

Based on this theorem, it is obvious how to give equivalent definitions of partial and total
degrees on all three layers of this theorem. For example, the $i^{\text{th}}$ partial degree of
a map $f\not\equiv0$ may be defined in the language of Statement\,\((i)\) by
\begin{equation}
\deg_i(f)+1\,:=\,\min\{\d\in\N\!\mit\!\underset{^i}\Delta^{\!\d}f\equiv0\}\,\ .
\end{equation}

\subsection{Lagrange Functions and Co"=Monofracts}

Next, we examine which functions $f\DP\Z_q\lmto\Z_r$ arise as polyfracts. For $d,q,r\in\N$
and $x\in\Z_q:=\Z/q\Z$, we set
\begin{equation}
\dbinom{d}{x}_{\!\!q,r}:=\,\sum_{0\leq\hat x\in x}(-1)^{\hat x}\dbinom{d}{\hat x}
+r\Z\,\in\,\Z_r\,\ ,\qquad\dbinom{d}{x}_{\!\!q}:=\,\dbinom{d}{x}_{\!\!q,0}\in\,\Z\,\ ,
\end{equation}
where only the at most $\cop{\tfrac{d}{q}}+1$ many summands with $\hat x\in
x\cap\{0,1,\dotsc,d\}\sb\Z$ may be nonzero. We call the map
\begin{equation}
\dbinom{d}{X}_{\!\!q,r}\DP\,\Z_q\lto\Z_r\,,\ x\lmto\dbinom{d}{x}_{\!\!q,r}\,\ ,\qquad
\dbinom{d}{X}_{\!\!q}:=\,\dbinom{d}{X}_{\!\!q,0}\,\ ,
\end{equation}
a \emph{co"=(mono)fract}. Most important will be the \emph{Lagrange Function}
\begin{equation}
\dbinom{0}{X}_{\!\!q,r}\!\DP\,x\lmto\dbinom{0}{x}_{\!\!q,r}=\,\begin{cases}
 1&\text{if $\,x=0$\,,}\\
 0&\text{if $\,x\neq0$\,,}
\end{cases}
\end{equation}
also written as
\begin{equation}
(1,0,0,\dotsc,0)\,\in\,\Z_r\!\´\n[q]=\Z_r^{\Z_q}\,\ ,
\end{equation}
and the shifted Lagrange Function
\begin{equation}
\dbinom{0}{X-x_0}_{\!\!q,r}\in\,\Z_r^{\Z_q}
\end{equation}
with center $x_0\in\Z_q$. Once more, Pascal's rule has nice consequences. Taking the
factor $(-1)^{\hat x}$ in our definition into account, we see that
\begin{equation}
\dbinom{d}{\ell+1}_{\!\!q}-\dbinom{d}{\ell}_{\!\!q}\,=\,\dbinom{d+1}{\ell+1}_{\!\!q}
\end{equation}
and
\begin{equation}\label{eq.DBo}
\Delta\dbinom{d}{X+x_0}_{\!\!q}\,=\,\dbinom{d+1}{X+x_0+1}_{\!\!q}\,\ .
\end{equation}
In the case of a prime power, $q=p^\ä$, $\ä\geq1$, we can use this to prove the
followinng\´:

\begin{Lemma}\label{lem.d}
Let $g\DP\Z_{p^\ä}\lto\Z$ be a function with vanishing value sum,
$\sum\limits_{\smash{x\in\Z_{p^\ä}}}g(x)=0$, then
$$
p\Div\Delta^{\!p^\ä-1}g\qquad\text{i.e.}\ \quad
p\Div\Delta^{\!p^\ä-1}g(x)\quad\text{for all $x\in\Z_{p^\ä}$.}
$$
\end{Lemma}

\begin{Beweis}
Since $\Delta$ is linear and goes well together with translations, it is enough to prove the
lemma for
\begin{sequation}
g=(g(0),g(1),\dotsc,g(p^\ä{-}1))\,:=\,(-1,0,0,\dotsc,0,1)\,\ .
\end{sequation}%
The translations of this function span the \(\Z\)"~submodule of functions $\Z_{p^\ä}\!\lto\Z$
with vanishing value sum. We also can write this function as
\begin{sequation}
g\,=\,\Delta\dbinom{0}{X}_{\!\!p^\ä}\,\ ,
\end{sequation}%
which makes our divisibility quite obvious in the point $x=0$:
\begin{sequation}
p\Div1+(-1)^{p^\ä}\,=\,\dbinom{p^\ä}{0}+(-1)^{p^\ä}\dbinom{p^\ä}{p^\ä}
 \,=\,\dbinom{p^\ä}{0+p^\ä}_{\!\!p^\ä}
\eqby{\eqref{eq.DBo}}\,\biggl.\Delta^{\!p^\ä}\dbinom{0}{X}_{\!\!p^\ä}\,\biggr|_{X=0}
 \!=\,\Delta^{\!p^\ä-1}g(0)\,\ .
\end{sequation}%
To prove the divisibility for $x\neq0$, let $\hat x\in x$ be the unique representative of $x$
with $0<\hx<p^\ä$, then
\begin{sequation}
p\Div\pm\dbinom{p^\ä}{\hx}\,=\,\dbinom{p^\ä}{x+p^\ä}_{\!\!p^\ä}
\eqby{\eqref{eq.DBo}}\,\Delta^{\!p^\ä}\dbinom{0}{X}_{\!\!p^\ä}\,\biggr|_{X=x}
\!=\,\Delta^{\!p^\ä-1}g(x)\,\ ,
\end{sequation}%
by Kummer's classical theorem about the multiplicity of primes in binomial coefficients.
\end{Beweis}

Since derivatives $\Delta f$ have vanishing value sum, we might apply this lemma repeatedly
to derivatives $\Delta f$ and would obtain the following\´:

\begin{Lemma}
For arbitrary functions $f\DP\Z_{p^\ä}\lto\Z$ and $\b\in\N$,
$$
p^\b\Div\Delta^{\!\b(p^\ä-1)+1}f\,\ .
$$
\end{Lemma}

\begin{Beweis}
The case $\b=0$ is trivial. To prove the induction step $\b\to\b{+}1$, we assume that the
formula holds for $\b$, so that the function
\begin{sequation}
g\,:=\,\dfrac{1}{p^\b}\´\Delta^{\!\b(p^\ä-1)+1}f
\end{sequation}%
is integer valued. Now, Lemma\,\ref{lem.d} applied to $g$ immediately yields the formula for
$\b+1$. We just have to use the trivial fact that any derivative $\Delta f$ has vanishing value
sum,
\begin{sequation}
\sum_{x\in\Z_q}\Delta f(x)\ =\ \sum_{x\in\Z_q}f(x+1)-f(x)
\ =\ \sum_{x\in\Z_q}f(x+1)-\sum_{x\in\Z_q}f(x)\ =\ 0\,\ .
\end{sequation}%
\end{Beweis}

However, the exponent $\b(p^\ä-1)+1$ of $\Delta$ in our lemma is not optimal for $\ä>1$.
Using only the case $\ä=1$ of this first approach, we improve it as follows:

\begin{Satz}\label{sz.D}
For arbitrary functions $f\DP\Z_{p^\ä}\lto\Z$ and $\b\in\N$,
$$
p^\b\Div\Delta^{\!(\b(p-1)+1)p^{\ä-1}}f\,\ .
$$
\end{Satz}

\begin{Beweis}
For $\gamma\in\N$ let
\begin{sequation}
p\nabla_{\!p^\gamma}\,:=\,\Delta^{\!p^\gamma}-\Delta_{p^\gamma}
\,=\,\sum_{i=1}^{p^\gamma-1}(-1)^{p^\gamma-i}\dbinom{p^\gamma}{i}\,T^i
\,\in\,\End(\Z_{p^\ä}^\Z)\,\ .
\end{sequation}%
Then $\nabla_{\!p^{\ä-1}}=\tfrac{1}{p}\,p\nabla_{\!p^{\ä-1}}$ lies in $\Z[T]$, as $p$ divides
$\tbinom{p^\gamma}{1}$, $\tbinom{p^\gamma}{2}$, \dots,
$\tbinom{p^\gamma}{p^\gamma-1}$. In particular, $\nabla_{\!p^{\ä-1}}$ commutes with
$\Delta_{p^{\ä-1}}\in\Z[T]$. Hence,
\begin{sequation}\label{eq.S}
\begin{split}
&\Delta^{\!(\b(p-1)+1)p^{\ä-1}}\,=\,\bigl(\Delta_{p^{\ä-1}}
     +p\nabla_{\!p^{\ä-1}}\bigr)^{\b(p-1)+1}\\[2pt]
&\ \ =\,\Delta_{p^{\ä-1}}^{\b(p-1)+1}+\,\bigl(p\nabla_{\!p^{\ä-1}}\bigr)^{\b(p-1)+1}
      +\,\sum_{i=1}^{\b(p-1)}\dbinom{\b(p-1)+1}{i}\Delta_{p^{\ä-1}}^{\b(p-1)+1-i}
      \bigl(p\nabla_{\!p^{\ä-1}}\bigr)^i\\[2pt]
&\ \ =\,\Delta_{p^{\ä-1}}^{\b(p-1)+1}
      +\,\bigl(p\nabla_{\!p^{\ä-1}}\bigr)^{\overbrace{^{\b(p-1)+1}}^{\geq\b+1}}
      +\,\sum_{\gamma=0}^{\b-1}\,\,\sum_{i=\gamma(p-1)+1}^{(\gamma+1)(p-1)}
 \,\,\,\overbrace{\!\!\!p^i}^{i\geq\gamma+1\!\!\!\!\!\!}
 \,\,\,\Delta_{p^{\ä-1}}^{\overbrace{^{\b(p-1)+1-i}}%
 ^{\geq(\b-(\gamma+1))(p-1)+1\!\!\!\!\!\!\!\!\!\!\!\!\!\!\!\!\!}}
 \,\,\,\dbinom{\b(p-1)+1}{i}\nabla_{\!p^{\ä-1}}^{\,i}\\
&\ \ =\,\Delta_{p^{\ä-1}}^{\b(p-1)+1}+\,p^{\b+1}S\,+\,\sum_{\gamma=0}^{\b-1}\,\,
p^{\gamma+1}\Delta_{p^{\ä-1}}^{(\b-(\gamma+1))(p-1)+1}S_\gamma\,\ ,
\end{split}
\end{sequation}%
with certain operators $S,S_0,S_1,\dotsc,S_{\b-1}\in\End(\Z_{p^\ä}^\Z)$. Now, if we apply
the operator $\Delta_{p^{\ä-1}}$ in this expression to a function $f\DP\Z_{p^\ä}\lto\Z$, it acts
within certain subdomains $D_j\sb\Z_{p^\ä}$. More precisely, we can partition the domain
$\Z_{p^\ä}$ into $p^{\ä-1}$ many \(p\)"~sets
\begin{sequation}
D_j\,:=\,\{j,j+p^{\ä-1},j+2p^{\ä-1},\dotsc,j+(p-1)p^{\ä-1}\}\qquad(j=0,1,\dotsc,p^{\ä-1}{-}1)\,\ .
\end{sequation}%
The restricted image $(\Delta_{p^{\ä-1}}f)|_{D_j}$ of $f$ under $\Delta_{p^{\ä-1}}$ depends
only on $f|_{D_j}$, and $\Delta_{p^{\ä-1}}$ acts on the restricted maps $f|_{D_j}$ in the
same way as $\Delta$ acts on functions $\Z_p\lto\Z$. Therefore, the last lemma, with
$\ä=1$, can be applied to the parts $f|_{D_j}$ of $f$. We obtain, for any $\gamma\in\N$ and
any $U\in\End(\Z_{p^\ä}^\Z)$,
\begin{sequation}
p^\gamma\Div\Delta^{\!\gamma(p-1)+1}_{p^{\ä-1}}f|_{D_j}
\ \ \text{,}\quad
p^\gamma\Div\Delta^{\!\gamma(p-1)+1}_{p^{\ä-1}}f
\ \quad\text{and}\quad
p^\gamma\Div\Delta^{\!\gamma(p-1)+1}_{p^{\ä-1}}U(f)\,\ .
\end{sequation}%
With Equation\,\eqref{eq.S} above, this yields
\begin{sequation}
p^\b\Div\Delta^{\!(\b(p-1)+1)p^{\ä-1}}f\,\ .
\end{sequation}%
\end{Beweis}

If we apply this Theorem to the \(\delta^{\text{th}}\) derivative $\binom{\delta}{X}_{\!\!p^\ä}$
of a co-monofract $\binom{0}{X-\delta}_{\!\!p^\ä}$, we obtain the following generalization of
Fleck's devisability relation \cite[Equation\,\((12)\)]{gr}\´:

\begin{Korollar}\label{cor.dpb}
If $\d\geq(\b(p-1)+1)p^{\ä-1}$ then
$$
p^\b\,\Div\,\dbinom{\d}{x}_{\!\!p^\ä}\ \quad\text{for all}\quad x\in\Z_{p^\ä}\,\ .
$$
\end{Korollar}

With the help of the last theorem, we obtain the following important result, which generalizes
\cite[Lemma\,1.5]{na}\´:

\begin{Satz}[Lagrange Polynomials]\label{sz.lagr}
For primes $p$, for $\ä,\b>0$ and $x_0\in\Z_{p^\ä}$
$$
\dbinom{0}{X-x_0}_{\!\!p^\ä\!\!,\,p^\b\!}=
\,\sum_{\d=0}^d\dbinom{\d}{\d-x_0}_{\!\!p^\ä\!\!,\,p^\b}\´\dbinom{X}{\d}
\,\,\in\,\Z_{p^\b}\dbinom{X}{\Z_{p^\ä}}\,\ ,
$$
with $d\,:=\,p^\ä{\!-}1+(\b{-}1)(p{-}1)p^{\ä-1}$ and $\d-x_0=\d-\widehat{x_0}+p^\ä\Z$. Here,
$\widehat{x_0}$ denotes the least nonnegative representative of $x_0$, and with that
notation
$$
\dbinom{0}{X-x_0}_{\!\!p^\ä\!\!,\,p\!}
=\,\,\sum_{\d=\widehat{x_0}}^{p^\ä-1}(-1)^{\d-\widehat{x_0}}\´
  \dbinom{\d}{\d-\widehat{x_0}}\dbinom{X}{\d}\,\,+\,p\´\Z^{\Z_{p^\ä}}
$$
and
$$
\dbinom{0}{X}_{\!\!p^\ä\!\!,\,p\!}=\,\,+\dbinom{X}{0}-\dbinom{X}{1}
+\dbinom{X}{2}\mp\dotsb+(-1)^{p^\ä-1}\´\dbinom{X}{p^\ä{-}1}\,\,+\,p\´\Z^{\Z_{p^\ä}}\,\ .
$$
\end{Satz}

\begin{Beweis}
Based on Theorem\,\ref{sz.D} this follows from our Taylor"~type Theorem\,\ref{sz.tay}
applied to the \(p^\ä\)"~periodic map
$\tbinom{0}{X{-}x_0}_{\!\!p^\ä\!\!,\,p^\b\!}\DP\Z\lto\Z_{p^\b\!}$, as, by
Equation\,\eqref{eq.DBo}\´, for all $\d\in\N$,
\begin{sequation}
\Bigl.\Delta^{\!\d}\dbinom{0}{X-x_0}_{\!\!p^\ä\!\!,\,p^\b\!}\,\Bigr|_{X=0}
\,=\,\dbinom{\d}{\d-x_0}_{\!\!p^\ä\!\!,\,p^\b\!}\,\ .
\end{sequation}%
For $\b=1$ this expression and our formula simplifies further, as
\begin{sequation}%
\dbinom{\d}{\d-x_0}_{\!\!p^\ä\!\!,\,p}\,=\begin{cases}
\tbinom{\d}{\d-\widehat{x_0}}+p\Z&\text{if \,$\d\geq\widehat{x_0}$\,,}\\[.2em]
0+p\Z&\text{if \,$\d<\widehat{x_0}$\,.}
\end{cases}
\end{sequation}
\end{Beweis}

As any map $\Z_{p^\ä}\!\lto\Z_{p^\b}$ is a linear combination of Lagrange Functions
$\tbinom{0}{X-x_0}_{\!\!p^\ä\!\!,\,p^\b\!}$, and Lagrange Functions are polyfractal, all maps
$\Z_{p^\ä}\!\lto\Z_{p^\b}$ are polyfractal. This was already obtained in
\cite[Corollary\,4.16]{hr} and in \cite[Lemma\,1]{wi}\´. We formulate it here for functions in
$n$ variables with arbitrary finite commutative codomains\´:

\begin{Korollar}\label{kor.papb}
Any \(p^\ä\)"~periodic map $\Z\lto\Z_{p^\b}$ can be described by a \(\Z_{p^\b}\)"~polyfract,
$$
\Z_{p^\b}\dbinom{X}{\Z_{p^\ä}}\,=\´\,\Z_{p^\b}^{\Z_{p^\ä}}\ .
$$
More generally, if $q_1,\dotsc,q_n\geq p$ are powers of a prime $p$ and $B$ is a finite
commutative \(p\)"~group, then any map $\Z_{q_1}\!\times\dotsm\times\Z_{q_n}\lto B$ is
polyfractal,
$$
B\dbinom{X_1\,,\,\dotsc,\,X_n}{\Z_{q_1}\!\times\dotsm\times\Z_{q_n}}
\,=\,B^{\Z_{q_1}\!\times\dotsm\times\Z_{q_n}}\ .
$$
\end{Korollar}\vsp

Actually, we can provide here the following stronger theorem\´:

\begin{Satz}[Interpolation Theorem]\label{sz.ip}
Let $p\in\N$ be prime, $\ä_1,\ä_2\dotsc,\ä_n,\b\geq1$ and
$$
f\DP\Z_{p^{\ä_1}}\!\times\Z_{p^{\ä_2}}\!\times\dotsm\times\Z_{p^{\ä_n}}\lto\,\Z_{p^\b}\,,\
x\lmto f(x)
$$
be any map. Then the polyfract
$$P\,=\,\sum_{\d\in[d]}P_\d
\dbinom{X_1,X_2,\dotsc,X_n}{\d_1\,,\,\d_2\,,\,\dotsc,\,\d_n}
\,\in\,\Z_{p^\b}\dbinom{X_1\,,\,X_2\,,\,\dotsc,\,X_n}{\Z_{p^{\ä_1}}\!\times\Z_{p^{\ä_2}}
 \!\times\dotsm\times\Z_{p^{\ä_n}}}\,\ ,
$$
with coefficients
$$
P_\d\,:=\!
\sum_{x\in\Z_{p^{\ä_1}}\times\dotsm\times\Z_{p^{\ä_n}}}
\dbinom{\d_1}{\d_1-x_1}_{\!\!p^{\ä_1}\!,\,p^\b}
\dbinom{\d_2}{\d_2-x_2}_{\!\!p^{\ä_2}\!,\,p^\b}\dotsm\
\dbinom{\d_n}{\d_n-x_n}_{\!\!p^{\ä_n}\!,\,p^\b}f(x)
$$
and summation range
$$
[d]\,:=\,\{0,1,\dotsc,d_1\}\times\dotsm\times\{0,1,\dotsc,d_n\}\ ,\quad d_j
\,:=\,p^{\ä_j}{\!-}1+(\b{-}1)(p{-}1)p^{\ä_j-1}\,\ ,
$$
interpolates the map $f$,
$$
P(x)\,=\,f(x)\ \quad\text{for all
\,$x\in\Z_{p^{\ä_1}}\!\times\Z_{p^{\ä_2}}\!\times\dotsm\times\Z_{p^{\ä_n}}$.}
$$
\end{Satz}\vsp

\begin{Beweis}
The polyfract,
\begin{sequation}
P(X)\,:=\!\sum_{x\in\Z_{p^{\ä_1}}\times\dotsm\times\Z_{p^{\ä_n}}}
\dbinom{0}{X_1-x_1}_{\!\!p^{\ä_1}\!,\,p^\b}\dbinom{0}{X_2-x_2}_{\!\!p^{\ä_2}\!,\,p^\b}\dotsm\
\dbinom{0}{X_n-x_n}_{\!\!p^{\ä_n}\!,\,p^\b}f(x)\,\ ,
\end{sequation}%
interpolates $f$, since, if we substitute any point
$\tx=(\tx_1,\dotsc,\tx_n)\in\Z_{p^{\ä_1}}\times\dotsm\times\Z_{p^{\ä_n}}$ into $X$, only the
summand with $x=\tx$ survives. The statement of the theorem follows, since, by
Teorem\,\ref{sz.lagr}\´, the polyfractal expansion of each of the \(n\)"~dimensional Lagrange
Functions
\begin{sequation}
\dbinom{0}{X_1-x_1}_{\!\!p^{\ä_1}\!,\,p^\b}\dbinom{0}{X_2-x_2}_{\!\!p^{\ä_2}\!,\,p^\b}\dotsm\
\dbinom{0}{X_n-x_n}_{\!\!p^{\ä_n}\!,\,p^\b}
\end{sequation}%
has \(\d\)"~coefficient
\begin{sequation}\label{eq.fac}
\dbinom{\d_1}{\d_1-x_1}_{\!\!p^{\ä_1}\!,\,p^\b}
\dbinom{\d_2}{\d_2-x_2}_{\!\!p^{\ä_2}\!,\,p^\b}\dotsm\
\dbinom{\d_n}{\d_n-x_n}_{\!\!p^{\ä_n}\!,\,p^\b}\,\ .
\end{sequation}%
\end{Beweis}

One can show (e.g.\ as in \cite[Corollary\,4.16]{hr}) that our Lagrange Polyfracts
$\tbinom{0}{X-x_0}_{\!\!p^\ä\!\!,\,p^\b}$ actually have degree
$d:=p^\ä{\!-}1+(\b{-}1)(p{-}1)p^{\ä-1}$\!, i.e.\
\begin{equation}
0\,\neq\,\dbinom{d}{d}_{\!\!p^\ä\!\!,\,p^\b}\!
=\,\dbinom{d}{d-x_0}_{\!\!p^\ä\!\!,\,p^\b}\,\ ,
\end{equation}
or
\begin{equation}\label{eq.inv}
p^{\b}\,\nDiv\,\dbinom{d}{d}_{\!\!p^\ä}\!
=\,\dbinom{d}{d-x_0}_{\!\!p^\ä}\,\ ,
\end{equation}
where the stated equality, with arbitrary $x_0\in\Z$, follows from the invariance of leading
coefficient of polynomials under the transformation $X\lmto X-x_0$. However, if we view the
Lagrange Functions
$\tbinom{0}{X}_{\!\!p^\ä\!\!,\,p^{\b'}\!}\in\Z_{p^{\b'}}\tbinom{X}{\Z_{p^\ä}}$, with increased
$\b'>\b\geq1$, modulo the smaller power $p^\b$\!, we rediscover
$\tbinom{0}{X}_{\!\!p^\ä\!\!,\,p^\b\!}$,
\begin{equation}
{\bmod}^{p^{\b'}}_{p^\b}\dbinom{0}{x}_{\!\!p^\ä\!\!,\,p^{\b'}\!}
\,=\,\dbinom{0}{x}_{\!\!p^\ä\!\!,\,p^\b\!}
\ \quad\text{for all $x\in\Z_{p^\ä}$,}
\end{equation}
where, for divisors $r$ of $r'\in\Z$,
\begin{equation}
{\bmod}^{r'}_r\DP\,\Z_{r'}\lto\Z_r\ ,\ \,x+r'\Z\lmto x+r\Z\,\ .
\end{equation}
Hence, the coefficients of
$\tbinom{0}{X}_{\!\!p^\ä\!\!,\,p^{\b'}\!}\in\Z_{p^{\b'}}\tbinom{X}{\Z_{p^\ä}}$ viewed modulo
$p^\b$ must equal those of
$\tbinom{0}{X}_{\!\!p^\ä\!\!,\,p^\b\!}\in\Z_{p^\b}\tbinom{X}{\Z_{p^\ä}}$, and the coefficients
$\tbinom{\d}{\d}_{\!\!p^\ä\!\!,\,p^{\b'}\!}$ of $\tbinom{0}{X}_{\!\!p^\ä\!\!,\,p^{\b'}\!}$ with
$\d>\deg\tbinom{0}{X}_{\!\!p^\ä\!\!,\,p^\b\!}=d$ must become zero. This holds also for the
shifted polyfracts $\tbinom{0}{X-x_0}_{\!\!p^\ä\!\!,\,p^\b\!}$, i.e.\
\begin{equation}\label{eq.dpb}
p^\b\,\Div\,\dbinom{\d}{\d-x_0}_{\!\!p^\ä}\ \quad\text{for}\quad\d
\,\geq\,p^{\ä-1}+\b(p{-}1)p^{\ä-1}\,\ ,
\end{equation}
as we have seen in Corollary\,\ref{cor.dpb} already. In particular, the leading coefficient of
$\tbinom{0}{X-x_0}_{\!\!p^\ä\!\!,\,p^\b\!}$ is close to be zero in $\Z_{p^\b}$,
\begin{equation}\label{eq.div}
p^{\b-1}\,\Div\,\dbinom{d}{d-x_0}_{\!\!p^\ä}\,\ .
\end{equation}
From the degree restriction for Laplace Functions it follows that (due to Teorem\,\ref{sz.1}\´)
$d:=p^\ä{\!-}1+(\b{-}1)(p{-}1)p^{\ä-1}$ is the best upper bound for the degree of all polyfracts
$P\in\Z_{p^\b}\tbinom{X}{\Z_{p^\ä}}$, an upper bound without any analog in the theory of
polynomials. In $n$ dimensions, for \(n\)"~dimensional polyfracts, we have the corresponding
restrictions of the partial degrees, but also the following restriction of the total degree, that is
stronger than the one that would arise from the restricted partial degrees alone\´:
\begin{Satz}\label{sz.md}
The polyfracts $P$ in
$\Z_{p^\b}\dbinom{X_1,\,X_2,\,\dotsc,\,X_n}{\Z_{p^{\ä_1}}\!\times\Z_{p^{\ä_2}}
\!\times\dotsm\times\Z_{p^{\ä_n}}}$ have bounded degree,
$$
\deg(P)\,\leq\,
\deg\Biggl(\dbinom{0}{X_1}_{\!\!p^{\ä_1}\!,\,p^\b\!}\dbinom{0}{X_2}_{\!\!p^{\ä_2}\!,\,p^\b\!}
 \!\dotsm\dbinom{0}{X_n}_{\!\!p^{\ä_n}\!,\,p^\b\!}\,\Biggr)
=\,\sum_{j=1}^np^{\ä_j}-n+(\b{-}1)(p{-}1)p^{\ä_{\text{max}}-1}\,,
$$
where $\ä_{\text{max}}:=\max\limits_{1\leq j\leq n}\ä_j$, and this bound is best possible.
\end{Satz}

\begin{Beweis}
Assume $\ä_1=\ä_{\text{max}}$, and let $d_1:=p^{\ä_1}{\!-}1+(\b{-}1)(p{-}1)p^{{\ä_1}-1}$\!.
Inside the polyfractal representation of the Lagrange Function
\begin{sequation}
\dbinom{0}{X_1}_{\!\!p^{\ä_1}\!,\,p^\b\!}\dbinom{0}{X_2}_{\!\!p^{\ä_2}\!,\,p^\b\!}
 \!\dotsm\dbinom{0}{X_n}_{\!\!p^{\ä_n}\!,\,p^\b\!}\,\ ,
\end{sequation}%
the coefficient of the monofract
\begin{sequation}
\dbinom{X_1\,,\,X_2\,,\,\dotsc\,,\,X_n\ }{d_1,p^{\ä_2}{-}1,\dotsc,p^{\ä_n}{-}1}
\end{sequation}%
is
\begin{sequation}\label{eq.coef}
\dbinom{d_1}{d_1}_{\!\!p^{\ä_1}\!,\,p^\b}\´
\dbinom{p^{\ä_2}-1}{p^{\ä_2}-1}_{\!\!p^{\ä_2}\!,\,p^\b}
\dotsm\dbinom{p^{\ä_n}-1}{p^{\ä_n}-1}_{\!\!p^{\ä_n}\!,\,p^\b}
\end{sequation}%
by Theorem\,\ref{sz.lagr}\´. Due to Relation\,\eqref{eq.inv}\´, this coefficient is not zero in
$\Z_{p^{\b}}$,
\begin{sequation}
p^{\b}\,\Div\hspace{-11pt}\diagup\,\,\dbinom{d_1}{d_1}_{\!\!p^{\ä_1}}\ ,\quad
p\,\Div\hspace{-11pt}\diagup\,\,\dbinom{p^{\ä_2}-1}{p^{\ä_2}-1}_{\!\!p^{\ä_2}}\ ,\ \dotsc\ ,\quad
p\,\Div\hspace{-11pt}\diagup\,\,\dbinom{p^{\ä_n}-1}{p^{\ä_n}-1}_{\!\!p^{\ä_n}}\,\ ,
\end{sequation}%
proving the best possible statement. However, this is the limit of how far we can go. If we
take any monofract of higher total degree in our Lagrange Polyfract, then its coefficient will
be zero by Corollary\,\ref{cor.dpb}\´. For example, if we increase $p^{\ä_2}-1$ only, then the
second nondivisiblity above will get lost, due to Corollary\,\ref{cor.dpb}\´. Therefore, as
$p^{\b-1}$ divides $\tbinom{d_1}{d_1}_{\!\!p^{\ä_1}}$, again by Corollary\,\ref{cor.dpb}\´, the
corresponding Coefficient\,\eqref{eq.coef} will become zero in $\Z_\b$. This proves the
equality in the theorem. The inequality follows from the fact that all (shifted) Lagrange
Polyfracts have the same degree, so that any function, as linear combination of Lagrange
Functions, has at most their degree.
\end{Beweis}

In principle, it also would be possible to work with \emph{fractal series} and to write
\begin{equation}
\dbinom{0}{X-x_0}_{\!\!p^\ä\!}
=\,\sum_{\d=0}^\infty\dbinom{\d}{\d-x_0}_{\!\!p^\ä\!}\´\dbinom{X}{\d}\,\ ,
\end{equation}
with uniform convergence (on the full domain $\Z$) in the \(p\)"~adic metric, as this equation
holds modulo any $p^\b$\!. Actually, if we plug in only nonnegative $x\in\Z$, then all
summands with $\d>x$ become zero, so that the equation holds without employing any
concept of convergence in these points. However, we just wanted to mention fractal series
here. They might be an elegant tool, but we do not use them in the present paper.

\subsection{Structure of Periodic Polyfracts in One Variable}
\label{sec.str}

In this subsection, we examine which maps between finite commutative groups are
\emph{polyfractal}, i.e.\ can be described by polyfracts. The following theorem of Hrykaj
\cite[Theorem\,1]{hr} gives first important information\´:

\begin{Satz}\label{sz.0}
Let $r\geq0$. For polyfracts $P=\sum\limits_{\d\in\N}
P_\d\tbinom{X}{\d}\in\Z_r\tbinom{X}{\Z}$ and numbers $q\geq1$ the following are
equivalent\´:
\begin{enumerate}[(i)]
\item \(P\DP\Z\lto\Z_r\)\, is \(q\)"~periodic.
\item \(\fa\d\in\N\DP\quad\dbinom{q}{1}P_{\d+1}+\dbinom{q}{2}P_{\d+2}
    +\dotsb+\dbinom{q}{q}P_{\d+q}\,=\,0\ \)\,.
\end{enumerate}\medspace
\end{Satz}

\begin{Beweis}
We have
\begin{sequation}
\Delta_q\ =\ T^q-\Id\ =\ (\Delta+\Id)^q-\Id\ =\ \sum_{j=1}^q\dbinom{q}{j}\Delta^j\,\ ,
\end{sequation}%
so that, for any $\d\in\N$,
\begin{equation}
\Delta_q\dbinom{X}{\d}\ =\ \sum_{j=1}^q\dbinom{q}{j}\Delta^j\dbinom{X}{\d}
\ =\ \sum_{j=1}^{\min\{q,\d\}}\dbinom{q}{j}\dbinom{X}{\d-j}
\end{equation}%
and
\begin{equation}%
\Delta_qP\ =\ \sum_{\d\in\N}\sum_{j=1}^{\min\{q,\d\}} P_\d\dbinom{q}{j}\dbinom{X}{\d-j}
\ =\ \sum_{\d\in\N}\,\Bigl[\,\sum_{j=1}^qP_{\d+j}\dbinom{q}{j}\Bigr]\dbinom{X}{\d}\,\ .
\end{equation}
Hence, by Theorem\,\ref{sz.1}\´,
\begin{sequation}
\Delta_qP\,=\,0
\quad\,\lEqi\quad\,\fa\d\in\N\DP\ \sum_{j=1}^qP_{\d+j}\dbinom{q}{j}\,=\,0\,\ .
\end{sequation}%
\end{Beweis}

In the special case $r=p$ prime and $q=p^\ä$, Part\,\((ii)\) simplifies, as
\begin{equation}
\dbinom{p^\ä}{1},\dbinom{p^\ä}{2},\dotsc,\dbinom{p^\ä}{p^\ä-1}\,\equiv\,0\pmod{p}\,\ ,
\end{equation}
and we obtain the following corollary, which also follows from Theorem\,\ref{sz.md} in one
dimension and a counting argument in connection with Corollary\,\ref{kor.papb}\´:

\begin{Korollar}\label{kor.p}
For primes $p$ and $\ä\geq0$
$$
\Z_p\dbinom{X}{\Z_{p^\ä}}\,=\,\bigl\{P\in\Z_p\tbinom{X}{\Z}\!\mit\!\deg(P)<p^\ä\bigr\}
$$
\end{Korollar}

In our next result, we show that Theorem\,\ref{sz.0} imposes strong restrictions on the
\emph{wavelength} (i.e.\ shortest periodicity) of polyfractal maps. For example, for polyfractal
maps into $\Z_r$, being \(12\)"~periodic without being \(4\)"~periodic or \(6\)"~periodic is
only possible if $6\div r$. We have the following general theorem\´:

\begin{Satz}\label{sz.q-per}
Let $r'$ be coprime to $r\geq2$\!. Any \(qr'\!\)"~periodic polyfract $P\DP\Z\lto\Z_r$ is already
\(q\)"~periodic,
$$
\Z_r\dbinom{X}{\Z_{qr'}}\,=\,\Z_{r}\dbinom{X}{\Z_{q}}\,\ .
$$
In particular, any \(r'\)"~periodic polyfract $P\DP\Z\lto\Z_r$ is constant,
$$
\Z_r\dbinom{X}{\Z_{r'}}\,=\,\Z_r\dbinom{X}{\Z_1}\,=\´\,\Z_r\,\ .
$$
\end{Satz}

\begin{Beweis}
If an \(r'\)"~periodic polyfract $P$ would have degree $\deg(P)>0$ then the equation in
Theorem\,\ref{sz.0}\´\((ii)\) with $\d:=\deg(P)-1\geq0$ would be violated. Hence, any
\(r'\)"~periodic polyfract $P$ is already constant. Therefore, in the general \(qr'\!\)"~periodic
case, $P\in\Z_r\tbinom{X}{\Z_{qp'}}$, each of the \(r'\!\)"~periodic polyfracts
\begin{sequation}
P(qX)\,,\ P(qX+1)\,,\ \dotsc\,,\ P(qX+q-1)
\end{sequation}%
is already constant, so that $P$ is \(q\)"~periodic.
\end{Beweis}

The theorem shows that over $\Z_{p^\b}$, $p$ prime, any periodic polyfract has a
\(p^\ä\)"~wave"|length, for some $\ä\geq0$. Conversely, we know from
Corollary\,\ref{kor.papb} that every map $\Z_{p^\ä}\lto\Z_{p^\b}$ is polyfractal. Since there
are $(p^\b)^{(p^\ä)}$ many such maps, there must also exist a periodic polyfract to any of
the $(p^\b)^{(p^\ä)}$ many choices for the first $p^\ä$ information coefficients
(Remark\,\ref{rem.inf}\´). Hence, every sequence $P_0,P_1,\dotsc,P_{p^\ä-1}\in\Z_{p^\b}$
can be extended to a finite sequence fulfilling condition \((ii)\) of Theorem\,\ref{sz.0}\´. This
extension is unique, by Theorem\,\ref{sz.1} and Remark\,\ref{rem.inf}\´. By
Corollary\,\ref{sz.md}\´, $(\b{-}1)(p{-}1)p^{\ä-1}$ additional periodicity coefficients suffice.

\subsection{Structure of Periodic Polyfracts in Several Variables}

In this section, we turn to polyfracts in several variables. We will see later that, even if we
only are interested in univariate functions, results about polyfracts in several variables are
helpful. We start with a reformulation of our last theorem for polyfracts in several variables\´:

\begin{Korollar}\label{cor.rem}
If $q_n$ is coprime to $\abs{B}$ then
$$
B\dbinom{X_1,\dotsc,X_{n-1},X_n}{\Z_{q_1}\!\times\dotsm\times\Z_{q_{n-1}}\!\times\Z_{q_n}}
\,=\,B\dbinom{X_1,\dotsc,X_{n-1}}{\Z_{q_1}\!\times\dotsm\times\Z_{q_{n-1}}}\,\ ,
$$
i.e., $X_n$ does not occur in any polyfract
 $P\in B\tbinom{X_1,\dotsc,X_{n-1},X_n}{\Z_{q_1}\!\times\dotsm\times\Z_{q_{n-1}}\!\times\Z_{q_n}}$.
\end{Korollar}

\begin{Beweis}
At any point $(x_1,x_2,\dotsc,x_{n-1})$ the polyfract
\begin{sequation}
P(x_1,x_2,\dotsc,x_{n-1},X_n)\,\in\,B\dbinom{X_n}{\Z_{q_n}}
\end{sequation}%
is constant, as $q_n$ is coprime to $r$. Hence, the polyfractal map $P$ does not depend on
$X_n$, and we may replace $X_n$ with $0$ without changing the map $P$\!. The obtained
polyfract $P|_{X_n=0}$ has to be equal to the original $P$\!, by the uniqueness of polyfractal
representations, so that there was no $X_n$ in $P$ in the first place.
\end{Beweis}

From this corollary we obtain the following basic insight\´:

\begin{Satz}[Decomposition Theorem]\label{sz.decomp}
Let $B_1$ and $B_2$ be finite commutative groups, and $q_1,q_2,\dotsc,q_n>1$. If
$q_1q_2\dotsm q_s$ is coprime to $\abs{B_2}$, and $q_{s+1}q_{s+2}\dotsm q_n$ is
coprime to $\abs{B_1}$, then
$$
(B_1\times B_2)\dbinom{X_1,\dotsc,X_n}{\Z_{q_1}\!\times\dotsm\times\Z_{q_n}}
\,=\,B_1\dbinom{X_1,\dotsc,X_s}{\Z_{q_1}\!\times\dotsm\times\Z_{q_s}}
\,\times\,B_2\dbinom{X_{s+1},\dotsc,X_n}{\Z_{q_{s+1}}\!\times\dotsm\times\Z_{q_n}}
$$
\end{Satz}

\begin{Beweis}
Apply Corollary\,\ref{cor.rem} to the two factors in the representation
\begin{sequation}
(B_1\times B_2)\dbinom{X_1,\dotsc,X_n}{\Z_{q_1}\!\times\dotsm\times\Z_{q_n}}
\,=\,B_1\dbinom{X_1,\dotsc,X_n}{\Z_{q_1}\!\times\dotsm\times\Z_{q_n}}
\,\times\,B_2\dbinom{X_1,\dotsc,X_n}{\Z_{q_1}\!\times\dotsm\times\Z_{q_n}}\,\ .
\end{sequation}%
\end{Beweis}

Of cause, the decomposition step in this theorem may be applied repeatedly. In combination
with Corollary\,\ref{kor.papb} this yields our main theorem, that we will be complemented by
Theorem\,\ref{sz.cp} in the next subsection\´:

\begin{Satz}[Classification of Polyfractal Maps]\label{sz.per}
Let $A$ and $B$ be finite commutative groups, and $p_1,\dotsc,p_t$ the prime divisors of
$\abs{A}\abs{B}$. For $i=1,2,\dotsc,t$, denote with $A_i$ resp.\,$B_i$ the (possibly trivial)
\(p_i\)"~primary component of $A$ resp.\,$B$. Assume
$A_i=\Z_{q_{i,1}}\!\times\dotsm\times\Z_{q_{i,n_i}}$, for certain powers
$q_{i,1},\dotsc,q_{i,n_i}$ of $p_i$, and set $n:=n_1+\dotsb+n_t$. Then
\begin{equation*}
\begin{split}
B\´\dbinom{X_1,\dotsc,X_n}{A}
\,=&\,\,B_1\dbinom{X_1,\dotsc,X_{n_1}}{A_1}\times\dotsb\times\,B_t
\dbinom{X_{n-n_t+1},\dotsc,X_n}{A_t}\\[.6em]
\,=\,B_1^{A_1}\!\times\dotsm\times B_t^{A_t}
\,:=&\,\,\bigl\{(a_1,\dotsc,a_t)\mto(f_1(a_1),\dotsc,f_t(a_t))\!\mit\!f_1\in B_1^{A_1}\!, \dotsc,
f_t\in B_t^{A_t}\bigr\}\ .
\end{split}
\end{equation*}
\end{Satz}

In this theorem, the domain $A$ is written as a direct product of $n$ cyclic groups of prime
power order,
\begin{equation}
A\,=\,A_1\!\times\dotsm\times A_t
\,=\,\bigl(\Z_{q_{1,1}}\!\times\dotsm\times\Z_{q_{1,n_1}}\bigr)\times\dotsm\times\bigl(\Z_{q_{t,1}}
 \!\times\dotsm\times\Z_{q_{t,n_t}}\bigr)\,\ .
\end{equation}
This means, any $x\in A$ is viewed as an \(n\)"~tuple $(x_1,\dotsc,x_n)$, and we have one
symbolic variable $X_j$ for each coordinate $x_j$ of $x$. If this variable belongs to the, say,
$i_j^{\text{\ th}}$ prime $p_{i_j}$, i.e.\ to one of the cyclic factors of the \(p_{i_j}\)"~primary
component $A_{i_j}$ of $A$, or, in other words, if $i_j$ is the smallest number with
\begin{equation}
j\,\leq\,\sum_{i=1}^{i_j}n_i\,\ ,
\end{equation}
then this variable only occurs in monomials with coefficients from the \(p_{i_j}\)"~group
\begin{equation}
\{0\}\times\dotsm\times\{0\}\times B_{i_j}\times \{0\}\times\dotsm\times\{0\}\,\ .
\end{equation}
We may say the variable $X_j$ belongs to $B_{i_j}$. Hence, it belongs to $B_{i_j}$ if and
only if it belongs to $A_{i_j}$, i.e. to one cyclic factor of $A_{i_j}$. In particular, only variables
that belong to the same prim may occur together in one monomial. This is a negative result,
showing that the trivial combinations $(a_1,\dotsc,a_t)\mto(f_1(a_1),\dotsc,f_t(a_t))$ of
maps $f_i\DP A_i\lto B_i$ are the only polyfractal maps $A\lto B$, if $n$ variables are used
in the described way.

\subsection{Splitting Variables}
\label{sec.spl}

In Theorem\,\ref{sz.per}\´, we have shown that the trivial combinations of maps $A_i\lto B_i$
are the only polyfractal maps $A\lto B$. At least, we have proven that this is true if we break
$A$ down into cyclic groups of prime power order and introduce one variable for each of
them. However, any other way to write $A$ as a product of cyclic groups yields the same
class of maps $A\lto B$ as polyfractal representable maps. This is because there is a way to
split cyclic factors, and the corresponding variables in a polyfract, without changing the
polyfractal map. The splitting procedure is described in the following theorem. It is
formulated for polyfracts in one variable $X$ only, but splitting a variable $X_j$ in a
multivariate polynomial works just the same\´:

\begin{Satz}[Splitting Variables]\label{sz.split}
Assume $r_1,r_2\geq2$ and let $q_1$ be coprime to $r_2$ and $q_2$ be coprime to $r_1$.
If
$$
P\,=\,(P_1,P_2)
\,=\,\biggl(\sum_{\d\in\N}c_{1,\d}\tbinom{X}{\d}\,,\,\sum_{\d\in\N}c_{2,\d}\tbinom{X}{\d}\biggr)
\,=\,\sum_{\d\in\N}(c_{1,\d}\´,c_{2,\d})\tbinom{X}{\d}
\,\in\,(\Z_{r_1}\!\times\Z_{r_2})\dbinom{X}{\Z}
$$
is \(q_1q_2\)"~periodic, then
$$
\check P\,:=\,\bigl(P_1(X_1),P_2(X_2)\bigr)
\,=\,\sum_{\d\in\N}(c_{1,\d},0)\tbinom{X_1}{\d}\,+\,\sum_{\d\in\N}(0,c_{2,\d})\tbinom{X_2}{\d}
\,\in\,(\Z_{r_1}\!\times\Z_{r_2})\dbinom{X_1,\,X_2}{\Z\times\Z}\,\ ,
$$
is \((q_1,q_2)\)"~periodic. We have a ring isomorphisms $P\,\lmto\,\check P$,
$$
(\Z_{r_1}\!\times\Z_{r_2})\dbinom{X}{\Z_{q_1q_2}}\lto\,
(\Z_{r_1}\!\times\Z_{r_2})\dbinom{X_1,\,X_2}{\Z_{q_1}\!\times\Z_{q_2}}\,,
$$
$$
(c_{1,\d},c_{2,\d})\tbinom{X}{\d}\,\lmto\,(c_{1,\d},0)\tbinom{X_1}{\d}\,+\,(0,c_{2,\d})\tbinom{X_2}{\d}\,,\
$$
with invers $
Q\,\lmto\,Q\rlap{\ensuremath{|}}\kern0em\lower2ex\hbox{\ensuremath{|}}\!\!%
_{\atop{X_1:=X^{\phantom{H}}}{X_2:=X^{\phantom{H}}}}\!\!. $\smallskip

With the splitting isomorphism $\vartheta\DP\Z_{q_1q_2}\!\lto\Z_{q_1}\!\times\Z_{q_2}$,
$x+q_1q_2\Z\,\lmto(x+q_1\Z\´,x+q_2\Z)$, we may also write
$$
\check P\,=\,P\nach\vartheta^{-1}\quad\text{and}\quad
 Q\rlap{\ensuremath{|}}\kern0em\lower2ex\hbox{\ensuremath{|}}\!\!%
_{\atop{X_1:=X^{\phantom{H}}}{X_2:=X^{\phantom{H}}}}\!\!\!=\,Q\nach\vartheta\,\ ,
$$
as equations between maps on $\Z_{q_1}\!\times\Z_{q_2}$ resp.\,$\Z_{q_1q_2}$.
\end{Satz}

\begin{Beweis}
By Theorem\,\ref{sz.q-per} and Theorem\,\ref{sz.decomp}\´, the range of the map
$P\lmto\check P$ is
\begin{sequation}
\Z_{r_1}\dbinom{X_1}{\Z_{q_1q_2}}\times\Z_{r_2}\dbinom{X_2}{\Z_{q_1q_2}}
\,=\,\Z_{r_1}\dbinom{X_1}{\Z_{q_1}}\times\Z_{r_2}\dbinom{X_2}{\Z_{q_2}}
\,=\,(\Z_{r_1}\!\times\Z_{r_2})\dbinom{X_1,\,X_2}{\Z_{q_1}\!\times\Z_{q_2}}\,\ ,
\end{sequation}%
i.e.\ $P\lmto\check P$ is surjective. As the ring homomorphism $Q\,\lmto\,
 Q\rlap{\ensuremath{|}}\kern0em\lower2ex\hbox{\ensuremath{|}}\!\!%
_{\atop{X_1:=X^{\phantom{H}}}{X_2:=X^{\phantom{H}}}}\!\!$ is a post"=invers (left inverse),
our map $P\lmto\check P$ is also injective and an isomorphism. Furthermore,
\begin{sequation}
Q\rlap{\ensuremath{|}}\kern0em\lower2ex
 \hbox{\ensuremath{|}}\!\!_{\atop{X_1:=X^{\phantom{H}}}{X_2:=X^{\phantom{H}}}}
 \!\!\!(x+q_1q_2\Z)
\,=\,Q\rlap{\ensuremath{|}}\kern0em\lower2ex
 \hbox{\ensuremath{|}}\!\!_{\atop{X_1:=X^{\phantom{H}}}{X_2:=X^{\phantom{H}}}}\!\!\!(x)
\,=\,Q(x,x)
\,=\,Q(x+q_1\Z\,,\,x+q_2\Z)
\,=\,Q\nach\vartheta\,(x+q_1q_2\Z)\,\ ,
\end{sequation}%
for any $x\in\Z$. This proves the last equation of the theorem. It also implies that we may
write any $P$ as $\check P\nach\vartheta$, so that the remaining second to the last
equation follows.
\end{Beweis}

In applications, cyclic domains might be most important, but in the study of polyfractal
representability, we prefer to work with their primary components. As we have seen, we
always can do this. Repeated application of Theorem\,\ref{sz.split} shows, we can restrict
ourselves to polyfractal maps
 $\Z_{q_1}\times\Z_{q_2}\times\dotsm\times\Z_{q_n}\lto B$ with prime powers $q_j$, and
with one variable $X_j$ for each cyclic factor $\Z_{q_j}$. Furthermore, the obtained finest
decomposition $\Z_{q_1}\times\Z_{q_2}\times\dotsm\times\Z_{q_n}$ of the domain $A$ is
essentially unique, as group theory tells us. In particular, we obtain the following complement
to our classification, Theorem\,\ref{sz.per}\´, which generalizes \cite[Theorem\,1]{fr2}\´:

\begin{Satz}[Classification of Polyfractal Maps]\label{sz.cp}
Let $A$ and $B$ be finite commutative groups, then
$$
\bigl\{f\in B^A\!\mit\!\text{$f$\ps is polyfractal}\,\bigr\}\,=\,\prod_{\text{\(p\) prime}}B_p^{A_p}\ ,
$$
where $A_p$ resp.\,$B_p$ is the \(p\)"~primary component of $A$ resp.\,$B$.
\end{Satz}

We conclude this section and the paper with the following illustrative example\´:

\begin{Beispiel}
Assume a map $f\DP\Z_{50}\lto\Z_{12}$ is polyfractal,
\begin{equation}
f\,\in\,\Z_{12}\dbinom{X}{\Z_{50}}\,\ .
\end{equation}
Then
\begin{equation}
\phi\nach f\nach\vartheta^{-1}\,\in\,(\Z_{2^2}\!\times\Z_{3^1}\!\times\Z_{5^0})
 \dbinom{X_1,\,X_2,\,X_3}{\Z_{2^1}\!\times\Z_{3^0}\!\times\Z_{5^2}}
\,=\,\Z_{4}^{\Z_{2}}\times\Z_{3}^{\Z_{1}}\times\Z_{1}^{\Z_{25}}
\,\cong\,\Z_{4}^{2}\times\Z_{3}\,\ ,
\end{equation}
with the Chinese Reminder Isomorphisms
\begin{equation}
\phi\DP\Z_{12}\lto\Z_{2^2}\!\times\Z_{3^1}\!\times\Z_{5^0}\qquad\text{and}\qquad
\vartheta\DP\Z_{50}\!\lto\Z_{2^1}\!\times\Z_{3^0}\!\times\Z_{5^2}
\end{equation}
from Section\,\ref{sec.mi} and Theorem\,\ref{sz.split}\´. Hence,
\begin{equation}
\phi\nach f\nach\vartheta^{-1}\,=\,(g,c,0)
\end{equation}
for a \(2\)"~periodic map $g\in\Z_4^{\Z_2}$\!, %$g\DP\Z\lto\Z_{4}$,
a constant (\(1\)"~periodic) map $c\in\Z_3^{\Z_1}=\Z_3^{\{0\}}$ %$c\DP\Z\lto\Z_{3}$
and the zero map $0\in\Z_1^{\Z_{25}}=\{0\}^{\Z_{25}}$\!. %$0\DP\Z\lto\{0\}$.
%In particular, $f$ is already \(2\)"~periodic. %, $f\DP\Z_{2}\lto\Z_{12}$.
Since, the splitting isomorphism
\begin{equation}
\phi\DP x+12\Z\lmto(x+4\Z\,,\,x+3\Z,\,x+1\Z)
\end{equation}
has the invers
\begin{equation}
\phi^{-1}\DP (a+4\Z\,,\,b+3\Z,\,0+1\Z)\lmto-3a+4b+12\Z\,\ ,
\end{equation}
it follows that
\begin{equation}
\begin{split}
f(x+50\Z)\,&=\,(\phi^{-1}\nach \phi\nach f\nach\vartheta^{-1}\nach\vartheta)(x+50\Z)\\
\,&=\,\phi^{-1}\bigl((\phi\nach f\nach\vartheta^{-1})(\vartheta(x+50\Z))\bigr)\\
\,&=\,\phi^{-1}\bigl((g,c,0))(x+2\Z\,,\,x+1\Z\,,\,x+25\Z)\bigr)\\
\,&=\,\phi^{-1}(g(x+2\Z),c(x+1\Z),0(x+25\Z))\\
%\,&=\,-3g(x+2\Z)+4c(x+1\Z)\\
\,&=\,-3\hat g(x+2\Z)+4\hat c\,\,+12\Z\,\ ,
\end{split}
\end{equation}
where $\hat g\in\{0,1,2,3\}^{\Z_2}\sb\Z^{\Z_2}$ and $\hat c\in\{0,1,2\}\sb\Z$ are the least
nonnegative representatives of $g\in\Z_4^{\Z_2}=\Z^{\Z_2}/4\Z^{\Z_2}$ and
$c(0+1\Z)\in\Z_3=\Z/3\Z$. Alternatively, we could have defined the multiplication with $-3$,
resp.\ $4$, as map from $\Z_4$, resp.\ $\Z_3$, into $\Z_{12}$. Then $f=-3g+4c$, if we view
the involved maps $f,g,c$ as maps on the common domain $\Z$.

We have seen that, the maps $f\DP\Z_{50}\lto\Z_{12}$ that can be written as polynomial
with rational coefficients are exactly the maps of the form
\begin{equation}
f(x+50\Z)\,=\,-3\hat g(x+2\Z)+4\hat c\,\,+12\Z\,\ ,
\end{equation}
with $\hat g\in\{0,1,2,3\}^{\Z_2}$ and $\hat c\in\{0,1,2\}$. In other words, the \(50\)"~periodic
functions $f\DP\Z\lto\Z_{12}$ that are polyfractal are precisely the \(2\)"~periodic functions
with $f(1)-f(0)\in\{\bar 0,\bar 3,\bar 6,\bar 9\}=3\Z_{12}\sb\Z_{12}$. These are the (quite
strong) necessary and sufficient restrictions for a map $f\DP\Z_{50}\lto\Z_{12}$ to be
polyfractal.

Moreover, if we actually want to write a polyfractal map $f\DP\Z_{50}\lto\Z_{12}$ as polyfract
in $\Z_{12}\tbinom{X}{\Z_{50}}$, we may use Theorem\,\ref{sz.ip} to calculate the involved
polyfract $g\in\Z_4\tbinom{X}{\Z_2}$ as
\begin{equation}
g\,=\,2a\binom{X}{2}-a\binom{X}{1}+b\ \quad\text{with $b=g(0)$ and $a=g(0)-g(1)$.}
\end{equation}
Using this (and that $-3\cdot2\equiv6\pmod{12}$), we obtain
\begin{equation}
f\,=\,6a\binom{X}{2}+3a\binom{X}{1}+d\ \quad\text{with $d=f(0)$ and $3a=f(1)-f(0)$.}
\end{equation}
\end{Beispiel}

%% ------------------------------------------------------------------------------
\vspace{1.2cm}

\noindent\textbf{Acknowledgement\´:}\\
We want to thank Mark B.\ Krugman and the continuing academic support team of the Xi’an
Jiaotong"=Liverpool University for their support.

%-------------------------------------------------------------------------------

\end{document}